\documentclass[11pt,singlespace,letter]{article}

\usepackage{times}
\usepackage[margin=1.1in]{geometry}
\usepackage{hyperref}
\usepackage{authblk}

\usepackage{amsmath,amssymb, amsfonts}
\usepackage{bm}
\usepackage{algorithm,algorithmic,bbm}

\newcommand{\Ex}{{\mathbb E}} 
\newcommand{\set}{\mathcal} 
\newcommand*{\QED}{\hfill\ensuremath{\blacksquare}}

\newtheorem{theorem}{Theorem}
\newtheorem{claim}[theorem]{Claim}

\newtheorem{lemma}[theorem]{Lemma}
\newtheorem{proposition}[theorem]{Proposition}

\usepackage[utf8]{inputenc} 
\usepackage[T1]{fontenc}    
\usepackage{hyperref}       
\usepackage{url}            
\usepackage{booktabs}       
\usepackage{amsfonts}       
\usepackage{nicefrac}       
\usepackage{microtype}      

\title{Optimal Cluster Recovery \\in the Labeled Stochastic Block Model}
\author[1]{Se-Young Yun}
\author[2]{Alexandre Proutiere}
\affil[1]{Los Alamos National Laboratory, Los Alamos, NM, USA}
\affil[2]{KTH, The Royal Institute of Technology, Stockholm, Sweden}

\begin{document}

\maketitle

\begin{abstract}
We consider the problem of community detection or clustering in the labeled Stochastic Block Model (LSBM) with a finite number $K$ of clusters of sizes linearly growing with the global population of items $n$. Every pair of items is labeled independently at random, and label $\ell$ appears with probability $p(i,j,\ell)$ between two items in clusters indexed by $i$ and $j$, respectively. The objective is to reconstruct the clusters from the observation of these random labels. 

Clustering under the SBM and their extensions has attracted much attention recently. Most existing work aimed at characterizing the set of parameters such that it is possible to infer clusters either positively correlated with the true clusters, or with a vanishing proportion of misclassified items, or exactly matching the true clusters. We find  the set of parameters such that there exists a clustering algorithm with at most $s$ misclassified items in average under the general LSBM and for any $s=o(n)$, which solves one open problem raised in \cite{abbe2015community}. We further develop an algorithm, based on simple spectral methods, that achieves this fundamental performance limit within $O(n \mbox{polylog}(n))$ computations and without the a-priori knowledge of the model parameters.  
\end{abstract}

\section{Introduction}

Community detection consists in extracting (a few) groups of similar items from a large global population, and has applications in a wide spectrum of disciplines including social sciences, biology, computer science, and statistical physics. The communities or clusters of items are inferred from the observed pair-wise similarities between items, which, most often, are represented by a graph whose vertices are items and edges are pairs of items known to share similar features. 

The stochastic block model (SBM), introduced three decades ago in \cite{holland1983}, constitutes a natural performance benchmark for community detection, and has been, since then, widely studied. In the SBM, the set of items $\mathcal{V}=\{1,\ldots,n\}$ are partitioned into $K$ non-overlapping clusters $\mathcal{V}_1,\ldots,\mathcal{V}_K$, that have to be recovered from an observed realization of a random graph. In the latter, an edge between two items belonging to clusters $\mathcal{V}_i$ and $\mathcal{V}_j$, respectively, is present with probability $p(i,j)$, independently of other edges. The analyses presented in this paper apply to the SBM, but also to the {\it labeled} stochastic block model (LSBM) \cite{heimlicher2012community}, a more general model to describe the similarities of items. There, the observation of the similarity between two items comes in the form of a {\it label} taken from a finite set ${\cal L}=\{0,1,\ldots,L\}$, and label $\ell$ is observed between two items in clusters $\mathcal{V}_i$ and $\mathcal{V}_j$, respectively, with probability $p(i,j,\ell)$, independently of other labels. The standard SBM can be seen as a particular instance of its labeled counterpart with two possible labels 0 and 1, and where the edges present (resp. absent) in the SBM correspond to item pairs with label 1 (resp. 0). The problem of cluster recovery under the LSBM consists in inferring the hidden partition $\mathcal{V}_1,\ldots,\mathcal{V}_K$ from the observation of the random labels on each pair of items.

Over the last few years, we have seen remarkable progresses for the problem of cluster recovery under the SBM (see \cite{GaoMZZ15} for an exhaustive literature review), highlighting its scientific relevance and richness. Most recent work on the SBM aimed at characterizing the set of parameters (i.e., the probabilities $p(i,j)$ that there exists an edge between nodes in clusters $i$ and $j$ for $1\le i,j\le K$) such that some qualitative recovery objectives can or cannot be met. For sparse scenarios where the average degree of items in the graph is $O(1)$, parameters under which it is possible to extract clusters positively correlated with the true clusters have been identified \cite{decelle2011, mossel2012stochastic, massoulie2013}. When the average degree of the graph is $\omega(1)$, one may predict the set of parameters allowing a cluster recovery with a vanishing (as $n$ grows large) proportion of misclassified items \cite{yun2014, mossel2014consistency}, but one may also characterize parameters for which an asymptotically exact cluster reconstruction can be achieved \cite{abbe2014exact,yun2014b,hajek2014achieving,mossel2014consistency,abbe2015community,abbe2015recovering,jog2015information}. 

In this paper, we address the finer and more challenging question of determining, under the general LSBM,  the minimal number of misclassified items given the parameters of the model. Specifically, for any given $s = o(n)$, our goal is to identify the set of parameters such that it is possible to devise a clustering algorithm with at most $s$ misclassified items. Of course, if we achieve this goal, we shall recover all the aforementioned results on the SBM. 

\medskip
\noindent
{\bf Main results.} We focus on the labeled SBM as described above, and where each item is assigned to cluster $\mathcal{V}_k$ with probability $\alpha_k>0$, independently of other items. We assume w.l.o.g. that $\alpha_1 \le \alpha_2 \le \dots \le \alpha_K$. We further assume that $\alpha=(\alpha_1,\ldots,\alpha_K)$ does not depend on the total population of items $n$. Conditionally on the assignment of items to clusters, the pair or edge $(v,w)\in \mathcal{V}^2$ has label $\ell\in {\cal L}=\{0,1,\dots,L\}$ with probability $p(i,j,\ell)$, when $v\in \mathcal{V}_i$ and $w\in \mathcal{V}_j$. W.l.o.g., 0 is the most frequent label, i.e., $0 = \arg \max_{\ell} \sum_{i=1}^{K}\sum_{j=1}^K \alpha_i \alpha_j p(i,j,\ell)$. Throughout the paper, we typically assume that $\bar{p}=o(1)$ and $\bar{p}n=\omega(1)$ where $\bar{p}=\max_{i,j,\ell\ge 1}p(i,j,\ell)$ denotes the maximum probability of observing a label different than 0. We shall explicitly state whether these assumption are made when deriving our results. In the standard SBM, the second assumption means that the average degree of the corresponding random graph is $\omega(1)$. This also means that we can hope to recover clusters with a vanishing proportion of misclassified items. We finally make the following assumption: there exist positive constants $\eta$ and $\varepsilon$ such that for every $i,j,k\in [K]=\{1,\ldots, K\}$, 
\begin{equation*}
\hbox{(A1)}\quad \forall \ell\in {\cal L}, \ \frac{p(i,j,\ell)}{p(i,k,\ell)} \le \eta\quad\hbox{and}\quad\hbox{ (A2)}\quad \frac{\sum_{k=1}^K\sum_{\ell=1}^L (p(i,k,\ell)-p(j,k,\ell))^2}{\bar{p}^2} \ge \varepsilon.
\end{equation*}
(A2) imposes a certain separation between the clusters. For example, in the standard SBM with two communities,  $p(1,1,1)=p(2,2,1)=\xi$, and $p(1,2,1)=\zeta$, (A2) is equivalent to $2(\xi-\zeta)^2/\xi^2\ge \epsilon$. In summary, the LSBM is parametrized by $\alpha$ and $p=(p(i,j,\ell))_{1\le i,j\le K, 0\le \ell\le L}$, and recall that $\alpha$ does not depend on $n$, whereas $p$ does.

For the above LSBM, we derive, for any arbitrary $s=o(n)$, a necessary condition under which there exists an algorithm inferring clusters with $s$ misclassified items. We further establish that under this condition, a simple extension of spectral algorithms extract communities with less than $s$ misclassified items. To formalize these results, we introduce the {\it divergence} of $(\alpha,p)$. We denote by $p(i)$ the $K\times (L+1)$ matrix whose element on the $j$-th row and the $(\ell+1)$-th column is $p(i,j,\ell)$ and denote by $p(i,j)\in [0,1]^{L+1}$ the vector describing the probability distribution of the label of a pair of items in $\mathcal{V}_i$ and $\mathcal{V}_j$, respectively. Let $\mathcal{P}^{K\times (L+1)}$
denote the set of $K \times (L+1)$ matrices such that each row represents
a probability distribution. The divergence $D(\alpha,p)$ of $(\alpha,p)$ is defined as follows: $D(\alpha,p) = \min_{i,j:i\neq j} D_{L+}(\alpha,p(i),p(j))$ with
\begin{align*}
D_{L+}(\alpha,p(i),p(j)) &= \min_{y\in \mathcal{P}^{K\times (L+1)}} \max\left\{ \sum_{k=1}^K \alpha_k KL(y(k), p(i,k)),\sum_{k=1}^K \alpha_k  KL(y(k), p(j,k))  \right\}
\end{align*}
where $KL$ denotes the Kullback-Leibler divergence between two label distributions, i.e.,\\ $KL(y(k),p(i,k)) = \sum_{\ell=0}^{L} y(k,\ell) \log\frac{y(k,\ell) }{ p(i,k,\ell)} $. Finally, we denote by $\varepsilon^\pi(n)$ the number of misclassified items under the clustering algorithm $\pi$, and by $\mathbb{E}[\varepsilon^\pi(n)]$ its expectation (with respect to the randomness in the LSBM and in the algorithm).

\medskip
\noindent
We first derive a tight lower bound on the average number of misclassified items when the latter is $o(n)$. Note that such a bound was unknown even for the SBM \cite{abbe2015community}.

\begin{theorem} Assume that (A1) and (A2) hold, and that $\bar{p}n=\omega(1)$. Let $s=o(n)$. If there exists a clustering algorithm $\pi$ misclassifying in average less than $s$ items asymptotically, i.e., $\lim\sup_{n\to\infty}{\mathbb{E}[\varepsilon^\pi(n)]\over s}\le 1$, then the parameters $(\alpha,p)$ of the LSBM satisfy: 
\begin{equation}\label{eq:loww}
\lim\inf_{n\to \infty} \frac{nD(\alpha,p)}{\log(n/s)} \ge 1.
\end{equation}
\label{thm:lower}
\end{theorem}

To state the corresponding positive result (i.e., the existence of an algorithm misclassifying only $s$ items), we make an additional assumption to avoid extremely sparse labels: (A3) there exists a constant $\kappa > 0$ such that
$np(j,i,\ell) \ge (n\bar{p})^\kappa$ for all $i,j$ and $\ell \ge 1$. 

\begin{theorem}\label{thm:algorithms} Assume that (A1), (A2), and (A3) hold, and that $\bar{p}=o(1)$, $\bar{p}n=\omega(1)$. Let $s=o(n)$. If the parameters $(\alpha,p)$ of the LSBM satisfy (\ref{eq:loww}), then the Spectral Partition ($SP$) algorithm presented in Section \ref{sec:algo} misclassifies at most $s$ items with high probability, i.e., $\lim_{n\to\infty}\mathbb{P}[\varepsilon^{SP}(n)\le s]=1$.
\end{theorem} 

These theorems indicate that under the LSBM with parameters satisfying (A1) and (A2), the number of misclassified items scales at least as $n\exp(-nD(\alpha, p)(1+o(1))$ under any clustering algorithm, irrespective of its complexity. They further establish that the Spectral Partition algorithm reaches this fundamental performance limit under the additional condition (A3). We note that the SP algorithm runs in polynomial time, i.e., it requires $O(n^2 \bar{p}\log(n))$ floating-point operations.


We further establish a necessary and sufficient condition on the parameters of the LSBM  for the existence of a clustering algorithm recovering the clusters exactly with high probability. Deriving such a condition was also open  \cite{abbe2015community}.

\begin{theorem}
Assume that (A1) and (A2) hold. If there exists a clustering algorithm that does not misclassify any item with high probability, then the parameters $(\alpha,p)$ of the LSBM satisfy: $\lim\inf_{n\to \infty} \frac{nD(\alpha,p)}{\log(n)} \ge 1$. If this condition holds, then under (A3), the SP algorithm recovers the clusters exactly with high probability.
\label{thm:exact}
\end{theorem}

The paper is organized as follows. Section 2 presents the related work and example of application of our results. In Section 3, we sketch the proof of Theorem 1, which leverages change-of-measure and coupling arguments. We present in Section 4 the Spectral Partition algorithm, and analyze its performance (we outline the proof of Theorem 2). All results are proved in details in the supplementary material.  

\section{Related Work and Applications}

\subsection{Related work}

Cluster recovery in the SBM has attracted a lot of attention recently. We summarize below existing results, and compare them to ours. Results are categorized depending on the targeted level of performance. First, we consider the notion of {\it detectability}, the lowest level of performance requiring that the extracted clusters are just positively correlated with the true clusters. Second, we look at {\it asymptotically accurate recovery}, stating that the proportion of misclassified items vanishes as $n$ grows large. Third, we present existing results regarding exact cluster recovery, which means that no item is misclassified. Finally, we report recent work whose objective, like ours, is to characterize the optimal cluster recovery rate.

\medskip
\noindent
{\bf Detectability.} Necessary and sufficient conditions for {\it detectability} have been
studied for the binary symmetric SBM (i.e., $L=1$, $K=2$, $\alpha_1=\alpha_2$,  $p(1,1,1)=p(2,2,1)=\xi$, and $p(1,2,1)=p(2,1,1)=\zeta$). In the sparse regime where $\xi,\zeta=o(1)$, and for the binary symmetric SBM, the main focus has been on identifying the phase transition threshold (a condition on $\xi$ and $\zeta$) for {\it detectability}: It was conjectured in \cite{decelle2011} that if $n(\xi-\zeta)<\sqrt{2n(\xi+\zeta)}$ (i.e., under the threshold), no algorithm can perform better than a simple random assignment of items to clusters, and above the threshold, clusters can
partially be recovered. The conjecture was recently proved in \cite{mossel2012stochastic} (necessary condition), and \cite{massoulie2013} (sufficient condition). The problem of detectability has been also recently studied in \cite{zhang2015community} for the asymmetric SBM with more than two clusters of possibly different sizes. Interestingly, it is shown that in most cases, the phase transition for detectability disappears. 

The present paper is not concerned with conditions for detectability. Indeed detectability means that only a strictly positive proportion of items can be correctly classified, whereas here, we impose that the proportion of misclassified items vanishes as $n$ grows large.

\medskip
\noindent 
{\bf Asymptotically accurate recovery.} A necessary and sufficient condition for asymptotically accurate recovery in the SBM (with any number of clusters of different but linearly increasing sizes) has been derived in \cite{yun2014} and \cite{mossel2014consistency}. Using our notion of divergence specialized to the SBM, this condition is $n D(\alpha,p)=\omega(1)$. Our results are more precise since the minimal achievable number of misclassified items is characterized, and apply to a broader setting since they are valid for the generic LSBM. 

\medskip
\noindent {\bf Asymptotically exact recovery.} Conditions for exact cluster recovery in the SBM have been also recently studied. \cite{abbe2014exact,mossel2014consistency,hajek2014achieving} provide a necessary and sufficient condition for asymptotically exact recovery in the binary symmetric SBM. For example, it is shown that when $\xi={a\log(n)\over n}$ and $\zeta={b\log(n)\over n}$ for $a>b$, clusters can be recovered exactly if and only if ${a+b\over 2}-\sqrt{ab}\ge 1$. In \cite{abbe2015community,abbe2015recovering}, the authors consider a more general SBM corresponding to our LSBM with $L=1$. They define CH-divergence as:
$$D_+ (\alpha,p(i),p(j)) =\frac{n}{\log (n)}\max_{\lambda\in[0,1]} \sum_{k=1}^K \alpha_k \left((1-\lambda) p(i,k,1)+\lambda p(j,k,1) - p(i,k,1)^{1-\lambda}p(j,k,1)^\lambda \right),$$
and show that $\min_{i\neq j} D_+ (\alpha, p(i),p(j)) >1$ is a necessary and
sufficient condition for asymptotically exact reconstruction. The following claim, proven in the supplementary material, relates $D_+$ to $D_{L+}$.
\begin{claim}
When $ \bar{p} = o(1)$, we have for all $i,j$: 
$$D_{L+} (\alpha, p(i),p(j))\stackrel{n\to\infty}{\sim}\max_{\lambda\in[0,1]} \sum_{\ell=1}^L \sum_{k=1}^K \alpha_k \left((1-\lambda) p(i,k,\ell)+\lambda p(j,k,\ell) - p(i,k,\ell)^{1-\lambda}p(j,k,\ell)^\lambda \right).$$ \label{thm:ch}
\end{claim}
Thus, the results in \cite{abbe2015community,abbe2015recovering} are obtained by applying Theorem~\ref{thm:exact} and Claim~\ref{thm:ch}.

In \cite{jog2015information}, the authors consider a symmetric labeled SBM where communities are balanced (i.e., $\alpha_k = \frac{1}{K}$ for all $k$) and where label probabilities are simply defined as $p(i,i,\ell) = p(\ell)$ for all $i$ and $p(i,j,\ell)=q(\ell)$ for all $i\neq j$. It is shown that $\frac{nI}{\log (n)} > 1$ is necessary and sufficient for asymptotically exact recovery, where $I = - \frac{2}{K} \log \left(\sum_{\ell=0}^L \sqrt{p(\ell)q(\ell)} \right)$. We can relate $I$ to $D(\alpha,p)$:
\begin{claim}
In the LSBM with $K$ clusters, if $\bar{p} = o(1)$, and for all $i,j,\ell$ such that $i\neq j$, $\alpha_i = \frac{1}{K}$, $p(i,i,\ell) = p(\ell)$, and $p(j,k,\ell)=q(\ell)$, we have: $D(\alpha,p) \stackrel{n\to\infty}{\sim}  - \frac{2}{K} \log \left(\sum_{\ell=0}^L \sqrt{p(\ell)q(\ell)} \right).$ \label{thm:di}
\end{claim}
Again from this claim, the results derived in \cite{jog2015information} are obtained by applying Theorem~\ref{thm:exact} and Claim~\ref{thm:di}.

\medskip
\noindent
{\bf Optimal recovery rate.} In \cite{deshpande2015asymptotic,mossel2015density}, the authors consider the binary SBM in the sparse regime where the average degree of items in the graph is $O(1)$, and identify the minimal number of misclassified items for very specific intra- and inter-cluster edge probabilities $\xi$ and $\zeta$. Again the sparse regime is out of the scope of the present paper. \cite{zhou2015,GaoMZZ15} are concerned with the general SBM corresponding to our LSBM with $L=1$, and with regimes where asympotically accurate recovery is possible. The authors first characterize the optimal recovery rate in a minimax framework. More precisely, they consider a (potentially large) set of possible parameters $(\alpha,p)$, and provide a lower bound on the expected number of misclassified items for the worst parameters in this set. Our lower bound (Theorem \ref{thm:lower}) is more precise as it is model-specific, i.e., we provide the minimal expected number of misclassified items for a given parameter $(\alpha,p)$ (and for a more general class of models). Then the authors propose a clustering algorithm, with time complexity $O(n^3\log(n))$, and achieving their minimax recovery rate. In comparison, our algorithm yields an optimal recovery rate $O(n^2 \bar{p} \log (n))$ for any given parameter $(\alpha,p)$, exhibits a lower running time, and applies to the generic LSBM.   

\subsection{Applications}

We provide here a few examples of application of our results, illustrating their versatility. In all examples, $f(n)$ is a function such that $f(n)=\omega(1)$, and $a,b$ are fixed real numbers such that $a>b$.

\medskip
\noindent{\bf The binary SBM.} Consider the binary SBM where the average item degree is $\Theta(f (n))$, and represented by a LSBM with parameters $L=1$, $K=2$, $\alpha=(\alpha_1,1-\alpha_1)$, $p(1,1,1)=p(2,2,1)=\frac{a f(n)}{n}$, and $p(1,2,1)=p(2,1,1)=\frac{b f(n)}{n}$. From Theorems 1 and 2, the optimal number of misclassified vertices scales as $n\exp(-g(\alpha_1,a,b)f(n)(1+o(1)))$ when $\alpha_1\le 1/2$ (w.l.o.g.) and where
$$g(\alpha_1,a,b):=\max_{\lambda\in [0,1]}(1-\alpha_1-\lambda+2 \alpha_1\lambda) a + (\alpha_1 + \lambda -2 \alpha\lambda) b - \alpha_1 a^\lambda b^{(1-\lambda)} - (1-\alpha_1) a^{(1-\lambda)}b^\lambda.$$
It can be easily checked that $g(\alpha_1,a,b) \ge g(1/2,a,b)=\frac{1}{2}(\sqrt{a}-\sqrt{b})^2 $ (letting $\lambda=\frac{1}{2}$). The worst case is hence obtained when the two clusters are of equal sizes. When $f(n)=\log(n)$, we also note that the condition for asymptotically exact recovery is $g(\alpha_1,a,b)\ge 1$.

\medskip
\noindent{\bf Recovering a single hidden community.} As in \cite{hajek2015information}, consider a random graph model with a hidden community consisting of $\alpha n$ vertices, edges between vertices belonging the hidden community are present with probability $\frac{a f (n)}{n}$, and edges between other pairs are present with probability $\frac{b f (n)}{n}$. This is modeled by a LSBM with parameters $K=2$, $L=1$, $\alpha_1 =\alpha$, $p(1,1,1)= \frac{a f (n)}{n}$, and $p(1,2,1)=p(2,1,1)=p(2,2,1)= \frac{b f (n)}{n}$. The minimal number of misclassified items when searching for the hidden community scales as $n\exp(-h(\alpha,a,b)f(n)(1+o(1)))$ where
$$
h(\alpha,a,b):=\alpha\left( a- (a-b)\frac{1+ \log(a-b) - \log(a\log(a/b))}{\log (a/b)} \right).
$$
When $f(n)=\log(n)$, the condition for asymptotically exact recovery of the hidden community is $h(\alpha,a,b)\ge 1$.

\medskip
\noindent{\bf Optimal sampling for community detection under the SBM.} Consider a dense binary symmetric SBM with intra- and inter-cluster edge probabilities $a$ and $b$. In practice, to recover the clusters, one might not be able to observe the entire random graph, but sample its vertex (here item) pairs as considered in \cite{yun2014}. Assume for instance that any pair of vertices is sampled with probability $\frac{\delta f(n)}{n}$ for some fixed $\delta>0$, independently of other pairs. We can model such scenario using a LSBM with three labels, namely $\times$, 0 and 1, corresponding to the absence of observation (the vertex pair is not sampled), the observation of the absence of an edge and of the presence of an edge, respectively, and with parameters for all $i,j\in\{ 1,2\}$, $p(i,j,\times)=1-\frac{\delta f(n)}{n}$, $p(1,1,1)=p(2,2,1)=a\frac{\delta f(n)}{n}$, and $p(1,2,1)=p(2,1,1)=b\frac{\delta f(n)}{n}$. The minimal number of misclassified vertices scales as $n\exp(-l(\delta,a,b)f(n)(1+o(1)))$ where
$
l:=\delta(1-\sqrt{ab}-\sqrt{(1-a)(1-b)}).$
When $f(n)=\log(n)$, the condition for asymptotically exact recovery is $l(\alpha,a_+,a_-,b_+,b_-)\ge 1$.

\medskip
\noindent{\bf Signed networks.} Signed networks \cite{leskovec2010signed,traag2009community} are used  in social sciences to model positive and negative interactions between individuals. These networks can be represented by a LSBM with three possible labels, namely 0, + and -, corresponding to the absence of interaction, positive and negative interaction, respectively. Consider such LSBM with parameters: $K=2$, $\alpha_1=\alpha_2$, $p(1,1,+)=p(2,2,+)=\frac{a_+ f(n)}{n}$, $p(1,1,-)=p(2,2,-)=\frac{a_- f(n)}{n}$, $p(1,2,+)=p(2,1,+)=\frac{b_+ f(n)}{n}$, and $p(1,2,-)=p(2,1,-)=\frac{b_- f(n)}{n}$, for some fixed $a_+,a_-,b_+,b_-$ such that $a_+>b_+$ and $a_- < b_-$. The minimal number of misclassified individuals here scales as $n\exp(-m(\alpha,a_+,a_-,b_+,b_-)f(n)(1+o(1)))$ where
$$  
m(\alpha,a_+,a_-,b_+,b_-):=\frac{1}{2}\left((\sqrt{a_+}-\sqrt{b_+})^2 + (\sqrt{a_-}-\sqrt{b_-})^2 \right).
$$
When $f(n)=\log(n)$, the condition for asymptotically exact recovery is $l(\alpha,a_+,a_-,b_+,b_-)\ge  1$.

\section{Fundamental Limits: Change of Measures through Coupling}

In this section, we explain the construction of the proof of Theorem \ref{thm:lower}. The latter relies on an appropriate {\it change-of-measure} argument, frequently used to identify upper performance bounds in online stochastic optimization problems \cite{lai1985}. In the following, we refer to $\Phi$, defined by parameters $(\alpha,p)$, as the true stochastic model under which all the observed random labels are generated, and denote by $\mathbb{P}_\Phi=\mathbb{P}$ (resp. $\mathbb{E}_\Phi[\cdot]=\mathbb{E}[\cdot]$) the corresponding probability measure (resp. expectation). In our change-of-measure argument, we construct a second stochastic model $\Psi$ (whose corresponding probability measure and expectation are $\mathbb{P}_\Psi$ and $\mathbb{E}_\Psi[\cdot]$, respectively). Using a change of measures from $\mathbb{P}_\Phi$ to $\mathbb{P}_\Psi$, we relate the expected number of misclassified items $\mathbb{E}_\Phi[\varepsilon^\pi(n)]$ under any clustering algorithm $\pi$ to the expected (w.r.t. $\mathbb{P}_\Psi$) log-likelihood ratio ${\cal Q}$ of the observed labels under $\mathbb{P}_\Phi$ and $\mathbb{P}_\Psi$. Specifically, we show that, roughly, $\log(n/\mathbb{E}_\Phi[\varepsilon^\pi(n)])$ must be smaller than $\mathbb{E}_\Psi[{\cal Q}]$ for $n$ large enough. 

\medskip
\noindent{\bf Construction of $\psi$.} Let $(i^\star, j^\star) = \arg\min_{i,j: i<j} D_{L+}(\alpha,p(i),p(j))$, and let $v^\star$ denote the smallest item index that belongs to cluster $i^\star$ or $j^\star$. If both $\mathcal{V}_{i^\star}$ and $\mathcal{V}_{j^\star}$ are empty, we define $v^\star =n$. Let $q\in {\cal P}^{K\times (L+1)}$ such that:
$
D(\alpha,p) = \sum_{k=1}^K \alpha_k KL(q(k),p(i^\star,k))  = \sum_{k=1}^K \alpha_k KL(q(k),p(j^\star,k)).
$
The existence of such $q$ is proved in Lemma~\ref{lem:dq} in the supplementary material. Now to define the stochastic model $\Psi$, we couple the generation of labels under $\Phi$ and $\Psi$ as follows.

1. We first generate the random clusters $\mathcal{V}_1,\ldots, \mathcal{V}_K$ under $\Phi$, and extract $i^\star$, $j^\star$, and $v^\star$. The clusters generated under $\Psi$ are the same as those generated under $\Phi$. For any $v\in \mathcal{V}$, we denote by $\sigma(v)$ the cluster of item $v$. 

2. For all pairs $(v,w)$ such that $v\neq v^\star$ and $w\neq v^\star$, the labels generated under $\Psi$ are the same as those generated under $\Phi$, i.e., the label $\ell$ is observed on the edge $(v,w)$ with probability $p(\sigma(v),\sigma(w),\ell)$.

3. Under $\Psi$, for any $v\neq v^\star$, the observed label on the edge $(v,v^\star)$ under $\Psi$ is $\ell$ with probability $q(\sigma(v),\ell)$.

Let $x_{v,w}$ denote the label observed for the pair $(v,w)$. We introduce $\mathcal{Q}$, the log-likelihood ratio of the observed labels under $\mathbb{P}_\Phi$ and $\mathbb{P}_{\Psi}$ as: 
\begin{equation}\label{eq:defL}
\mathcal{Q} = \sum_{v=1}^{v^{\star}-1} \log \frac{q(\sigma(v),x_{v^{\star},v})}{p(\sigma(v^{\star}),\sigma(v),x_{v^{\star},v})}+\sum_{v=v^{\star}+1}^n \log \frac{q(\sigma(v),x_{v^{\star},v})}{p(\sigma(v^{\star}),\sigma(v),x_{v^{\star},v})}.
\end{equation}

Let $\pi$ be a clustering algorithm with output $(\hat{\mathcal V}_k)_{1\le k \le K}$, and let ${\cal E} = \bigcup_{1\le k \le K}\hat{\mathcal V}_k \setminus \mathcal{V}_k$ be the set of misclassified items under $\pi$. Note that in general in our analysis, we always assume without loss of generality that $| \bigcup_{1\le k \le K}\hat{\mathcal V}_k \setminus \mathcal{V}_k | \le |\bigcup_{1\le  k \le K} \hat{\mathcal V}_{\gamma(k)} \setminus \mathcal{V}_{k}|$ for any permutation $\gamma$, so that the set of misclassified items is indeed ${\cal E}$. By definition, $\varepsilon^\pi (n) = |{\cal E}|$. Since under $\Phi$, items are interchangeable (remember that items are assigned to the various clusters in an i.i.d. manner), we have: $n \mathbb{P}_{\Phi}\{v\in {\cal E} \}= \mathbb{E}_{\Phi}[\varepsilon^\pi (n)] = \mathbb{E}[\varepsilon^\pi (n)]$.

Next, we establish a relationship between $\Ex [\varepsilon^\pi (n) ]$ and the distribution of $\mathcal{Q}$ under $\mathbb{P}_\Psi$. For any function $f(n)$, we can prove that:
$
\mathbb{P}_{\Psi} \{ \mathcal{Q} \le  f(n) \} \le
\exp(f(n))\frac{\Ex_{\Phi}[\varepsilon^\pi(n) ]}{(\alpha_{i^\star}+\alpha_{j^\star}) n} +
\frac{\alpha_{j^\star}}{\alpha_{i^\star}+\alpha_{j^\star}}.\label{eq:2v1} 
$
Using this result with $f(n) = \log\left(n/\Ex_{\Phi}[\varepsilon^\pi(n) ] \right) - \log(2/\alpha_{i^\star} )$, and Chebyshev's inequality, we deduce that:
$
\log\left(n/\Ex_{\Phi}[\varepsilon^\pi(n) ] \right) - \log(2/\alpha_{i^\star} ) \le \mathbb{E}_{\Psi}[\mathcal{Q}] + \sqrt{\frac{4}{\alpha_{i^\star}} \mathbb{E}_{\Psi}[(\mathcal{Q}-\mathbb{E}_{\Psi}[\mathcal{Q}])^2]},
$
and thus, a necessary condition for $\mathbb{E}[\varepsilon^\pi(n)]\le s$ is:
\begin{equation}
\log\left(n/s \right) -
\log(2/\alpha_{i^\star} ) \le \mathbb{E}_{\Psi}[\mathcal{Q}] + \sqrt{\frac{4}{\alpha_{i^\star}} \mathbb{E}_{\Psi}[(\mathcal{Q}-\mathbb{E}_{\Psi}[\mathcal{Q}])^2]}. \label{eq:fnl0}
\end{equation}

\medskip
\noindent{\bf Analysis of ${\cal Q}$.}
In view of (\ref{eq:fnl0}), we can obtain a necessary condition for $\mathbb{E}[\varepsilon^\pi(n)]\le s$ if we evaluate $\mathbb{E}_{\Psi}[\mathcal{Q}]$ and $\mathbb{E}_{\Psi}[(\mathcal{Q}-\mathbb{E}_{\Psi}[\mathcal{Q}])^2]$. To evaluate $\mathbb{E}_{\Psi}[\mathcal{Q}]$, we can first prove that $v^\star\le \log(n)^2$ with high probability. From this, we can approximate $\mathbb{E}_{\Psi}[\mathcal{Q}]$ by $\mathbb{E}_{\Psi}[\sum_{v=v^{\star}+1}^n \log \frac{q(\sigma(v),x_{v^{\star},v})}{p(\sigma(v^{\star}),\sigma(v),x_{v^{\star},v})}]$, which is itself well-approximated by $nD(\alpha,p)$. More formally, we can show that: 
\begin{eqnarray}
\mathbb{E}_{\Psi}[\mathcal{Q}] &\le& \left(n+2\log(\eta)\log (n)^2 \right)D(\alpha,p) + \frac{\log\eta}{n^3}.\label{eq:90}
\end{eqnarray}
Similarly, we prove that $\mathbb{E}_{\Psi}[(\mathcal{Q}-\mathbb{E}_{\Psi}[\mathcal{Q}])^2] = O(n \bar{p})$, which in view of Lemma \ref{lem:dpq} (refer to the supplementary material) and assumption (A2), implies that: $\mathbb{E}_{\Psi}[(\mathcal{Q}-\mathbb{E}_{\Psi}[\mathcal{Q}])^2] = o(n D(\alpha,p))$.

We complete the proof of Theorem \ref{thm:lower} by putting the above arguments together: From \eqref{eq:fnl0}, \eqref{eq:90} and the above analysis of ${\cal Q}$, when the expected number of misclassified items is less than $s$ (i.e., $\Ex[\varepsilon^\pi(n) ] \le s$), we must have: $\liminf_{n\to \infty} \frac{nD(\alpha,p)}{\log\left(n/s \right)}\ge 1.$

\section{The Spectral Partition Algorithm and its Optimality}\label{sec:algo}

In this section, we sketch the proof of Theorem \ref{thm:algorithms}. To this aim, we present the Spectral Partition (SP) algorithm and analyze its performance. The SP algorithm consists in two parts, and its detailed pseudo-code is presented at the beginning of the supplementary document (see Algorithm \ref{alg:partition}).

The first part of the algorithm can be interpreted as an initialization for its second part, and consists in applying a spectral decomposition of a $n\times n$ random matrix $A$ constructed from the observed labels. More precisely, $A=\sum_{\ell =1}^L w_\ell A^\ell$, where $A^\ell$ is the binary matrix identifying the item pairs with observed label $\ell$, i.e., for all $v,w\in {\cal V}$, $A_{vw}^\ell =1$ if and only if $(v,w)$ has label $\ell$. The weight $w_\ell$ for label $\ell\in \{1,\ldots,L\}$ is generated uniformly at random in $[0,1]$, independently of other weights. From the spectral decomposition of $A$, we estimate the number of communities and provide asymptotically accurate estimates $S_1,\ldots,S_K$ of the hidden clusters asymptotically accurately, i.e., we show that when $n\bar{p}=\omega(1)$, with high probability, $\hat{K}=K$ and there exists a permutation $\gamma$ of $\{1,\ldots,K\}$ such that $\frac{1}{n}\left|\cup_{k=1}^K \mathcal{V}_k\setminus S_{\gamma(k)} \right| = O\left(\frac{\log(n\bar{p})^2}{n\bar{p}}\right)$. This first part of the SP algorithm is adapted from algorithms proposed for the standard SBM in \cite{coja2010, yun2014} to handle the additional labels in the model without the knowledge of the number $K$ of clusters. 

The second part is novel, and is critical to ensure the optimality of the SP algorithm. It consists in first constructing an estimate $\hat{p}$ of the true parameters $p$ of the model from the matrices $(A^{\ell})_{1\le \ell \le L}$ and the estimated clusters $S_1,\ldots,S_K$ provided in the first part of SP. We expect $p$ to be well estimated since $S_1,\ldots,S_K$ are asymptotically accurate. Then our cluster estimates are iteratively improved. We run $\lfloor\log(n)\rfloor$ iterations. Let $S_1^{(t)},\ldots,S_K^{(t)}$ denote the clusters estimated after the $t$-th iteration, initialized with $(S_1^{(0)},\ldots,S_K^{(0)})=(S_1,\ldots,S_K)$. The improved clusters $S_1^{(t+1)},\ldots,S_K^{(t+1)}$ are obtained by assigning each item $v\in \mathcal{V}$ to the cluster maximizing a log-likelihood formed from $\hat{p}$, $S_1^{(t)},\ldots,S_K^{(t)}$, and the observations $(A^{\ell})_{1\le \ell \le L}$: $v$ is assigned to $S_{k^\star}^{(t+1)}$ where $k^{\star} = \arg \max_{k} \{\sum_{i=1}^K \sum_{w \in S^{(t-1)}_i}\sum_{\ell =0}^L A_{vw}^\ell \log \hat{p}(k,i,\ell)  \} $.


\medskip
\noindent{\bf Part 1: Spectral Decomposition.}
The spectral decomposition is described in Lines 1 to 4 in Algorithm \ref{alg:partition}. As usual in spectral methods, the matrix $A$ is first trimmed (to remove lines and columns corresponding to items with too many observed labels -- as they would perturb the spectral analysis). To this aim, we estimate the average number of labels per item, and use this estimate, denoted by $\tilde{p}$ in Algorithm \ref{alg:partition}, as a reference for the trimming process. $\Gamma$ and $A_\Gamma$ denote the set of remaining items after trimming, and the corresponding trimmed matrix, respectively. 

If the number of clusters $K$ is known and if we do not account for time complexity, the two step algorithm in \cite{coja2010} can extract the clusters from $A_\Gamma$: first the optimal rank-$K$ approximation $A^{(K)}$ of $A_{\Gamma}$ is derived using the SVD; then, one applies the $k$-mean algorithm to the columns of $A^{(K)}$ to reconstruct the clusters. The number of misclassified items after this two step algorithm is obtained as follows.
Let $M^\ell=\mathbb{E}[A_\Gamma^\ell]$, and $M=\sum_{\ell=1}^Lw_\ell M^\ell$ (using the same weights as those defining $A$). Then, $M$ is of rank $K$.
If $v$ and $w$ are in the same cluster, $M_v= M_w$ and if $v$ and $w$ do not belong to the same cluster, from (A2), we must have with high probability: $\|M_v-M_w\|_2=\Omega(\bar{p}\sqrt{n})$. Thus, the $k$-mean algorithm misclassifies $v$ only if $\| A^{(K)}_v-M_v\|_2 = \Omega(\bar{p}\sqrt{n})$. By leveraging elements of random graph and random matrix theories, we can establish that $\sum_{v}\| A^{(k)}_v-M_v \|_2^2 = \| A^{(k)}-M \|_F^2=O(n\bar{p})$ with high probability. Hence the algorithm misclassifies $O(1/\bar{p})$ items with high probability. 

Here the number of clusters $K$ is not given a-priori. In this scenario, Algorithm~\ref{alg:specg} estimates the rank of $M$ using a singular value thresholding procedure. To reduce the complexity of the algorithm, the singular values and singular vectors are obtained using the iterative power method instead of a direct SVD. It is known from \cite{halko2011finding} that with $\Theta\left(\log (n)\right)$ iterations, the iterative power method find singular values and the rank-$K$ approximation very accurately. Hence, when $n\bar{p}=\omega(1)$, we can easily estimate the rank of $M$ by looking at the number of singular values above the threshold $\sqrt{n\tilde{p}}\log(n\tilde{p})$, since we know from random matrix theory that the $(K+1)$-th singular value of $A_{\Gamma}$ is much less than $\sqrt{n\tilde{p}}\log(n\tilde{p})$ with high probability. In the pseudo-code of Algorithm~\ref{alg:specg}, the estimated rank of $M$ is denoted by $\tilde{K}$. 

The rank-$\tilde{K}$ approximation of $A_\Gamma$ obtained by the iterative power method is $\hat{A}=\hat{U}\hat{V}= \hat{U}\hat{U}^\top A_\Gamma$. From the columns of $\hat{A}$, we can estimate the number of clusters and classify items. Almost every column of $\hat{A}$ is located around the corresponding column of $M$ within a distance $\frac{1}{2}\sqrt{\frac{n\tilde{p}^2}{\log (n\tilde{p})}}$, since $\sum_{v}\| \hat{A}_v-M_v \|_2^2 = \| \hat{A}-M \|_F^2=O(n\bar{p}\log(n\bar{p})^2)$ with high probability (we rigorously analyze this distance in the supplementary material Section~\ref{sec:part-1-sp}). From this observation, the columns can be categorised into $K$ groups.  To find these groups, we randomly pick $\log (n)$ reference columns and for each reference column, search all columns within distance $\sqrt{\frac{n\tilde{p}^2}{\log (n\tilde{p})}}$. Then, with high probability, each cluster has at least one reference column and each reference column can find most of its cluster members. Finally, the $K$ groups are identified using the reference columns. To this aim, we compute the distance of $n\log (n)$ column pairs $\hat{A}_v$, $\hat{A}_w$. Observe that $\|\hat{A}_v - \hat{A}_w \|_2 = \|\hat{V}_v - \hat{V}_w \|_2 $ for any $u,v\in \Gamma$, since the columns of $\hat{U}$ are orthonormal. Now $\hat{V}_v$ is of dimension $\tilde{K}$, and hence we can identify the groups using  $O(n\tilde{K}\log(n))$ operations.

\begin{theorem} Assume that (A1) and (A2) hold, and that $n\bar{p}=\omega(1)$. 
After Step 4 (spectral decomposition) in the SP algorithm, with high probability, $\hat{K}=K$ and there exists a permutation $\gamma$ of $\{1,\ldots,K\}$ such that:
$\left|\cup_{k=1}^K \mathcal{V}_k\setminus S_{\gamma(k)} \right| = O\left(\frac{\log(n\bar{p})^2}{\bar{p}}\right).$
\label{thm:spec}
\end{theorem}

\medskip
\noindent{\bf Part 2: Successive clusters improvements.} 
Part 2 of the SP algorithm is described in Lines 5 and 6 in Algorithm \ref{alg:partition}. To analyze the performance of each improvement iteration, we introduce the set of items $H$ as the largest subset of $\mathcal{V}$ such that for all $v\in H$:
(H1) $e(v,\mathcal{V}) \le 10 \eta n\bar{p}L$; (H2) when $v\in \mathcal{V}_k$, $ \sum_{i =1}^K \sum_{\ell =0}^L e (v,\mathcal{V}_i ,\ell ) \log\frac{p (k,i,\ell)}{p (j,i,\ell)}  \ge \frac{n\bar{p}}{\log (n\bar{p})^4}$ for all $j \neq k$; (H3) $e(v,\mathcal{V}\setminus H) \le 2 \log (n\bar{p})^2$, where for any $S\subset \mathcal{V}$ and $\ell$, $e(v,S,\ell)=\sum_{w\in S}A_{vw}^\ell$, and $e(v,S)=\sum_{\ell=1}^Le(v,S,\ell)$. Condition (H1) means that there are not too many observed labels $\ell\ge 1$ on pairs including $v$, (H2) means that an item $v\in \mathcal{V}_k$ must be classified to $\mathcal{V}_k$ when considering the log-likelihood, and (H3) states that $v$ does not share too many labels with items outside $H$.  

We then prove that $|\mathcal{V}\setminus H| \le s$ with high probability when $n D(\alpha,p)  - \frac{n\bar{p}}{\log (n\bar{p})^3} \ge  \log(n/s) + \sqrt{\log(n/s)}$. This is mainly done using concentration arguments to relate the quantity \\
$\sum_{i =1}^K \sum_{\ell =0}^L e (v,\mathcal{V}_i ,\ell ) \log\frac{p (k,i,\ell)}{p (j,i,\ell)}$ involved in (H2) to $nD(\alpha,p)$.     

Finally, we establish that if the clusters provided after the first part of the SP algorithm are asymptotically accurate, then after $\log(n)$ improvement iterations, there is no misclassified items in $H$. To that aim, we denote by $\mathcal{E}^{(t)}$ the set of misclassified items after the $t$-th iteration, and show that with high probability, for all $t$, $ \frac{ | \mathcal{E}^{(t+1)}\cap H|}{ |\mathcal{E}^{(t)}\cap H|} \le \frac{1}{\sqrt{n\bar{p}}}$. This completes the proof of Theorem \ref{thm:algorithms}, since after $\log(n)$ iterations, the only misclassified items are those in $\mathcal{V}\setminus H$.

\newpage

\bibliographystyle{plain}
\bibliography{reference}

\newpage

\appendix

\section{The SP Algorithm}

In this section, we present the Spectral Partition (SP) algorithm. The main pseudo-code of SP is presented in Algorithm \ref{alg:partition}. The SP algorithm consists in two parts. In the first part, corresponding to Lines 1-4 in the pseudo-code, we apply a spectral decomposition of the matrix $A=\sum_{\ell=1}^Lw_\ell A^\ell$ constructed from the observed labels. This matrix is first trimmed, and then treated by applying the spectral decomposition algorithm, whose pseudo-code is presented in Algorithm \ref{alg:specg}. The second part of the SP algorithm, corresponding to Lines 5 and 6 in Algorithm \ref{alg:partition}, consists in improving the clusters initially identified in the first step.

\begin{algorithm}[htb]
   \caption{Spectral Partition}
   \label{alg:partition}
\begin{algorithmic}
\STATE {\bfseries Input:} Observation matrices $A^{\ell}$ for every label $\ell$ ($A_{uv}^\ell = 1$ if $\ell$ is observed between $u$ and $v$).
\STATE {\bf 1. Estimated average degree.} $\tilde{p} \leftarrow \frac{\sum_{\ell=1}^{L}\sum_{u,v}A_{uv}^\ell}{n(n-1)}$
\STATE {\bf 2. Random Weights.} $A \leftarrow \sum_{\ell=1}^{L}w_\ell A^{\ell}$ where the weights $w_\ell$'s are i.i.d and uniformly distributed on $[0,1]$.
\STATE {\bf 3. Trimming.} Construct $A_{\Gamma}=(A_{vw})_{v,w\in \Gamma}$, where $\Gamma$ is the set of nodes obtained after removing $\lfloor n \exp(-n\tilde{p})\rfloor$ nodes having the largest $\sum_{\ell}\sum_{w \in V} A_{vw}^\ell$.
\STATE {\bf 4. Spectral Decomposition.} Run Algorithm~\ref{alg:specg}  with input $A_{\Gamma}, \tilde{p}$, and output $(S_k)_{k=1,\ldots,\hat{K}}$.
\STATE {\bf 5. Estimated parameters.} $\hat{p}(i,j,\ell) \leftarrow \frac{\sum_{u\in S_i}\sum_{v\in S_j}A^{\ell}_{uv}}{|S_i||S_j|}$ for all $1\le i,j \le \hat{K}$ and $0\le\ell\le L$.
\STATE {\bf 6. Improvement.}
   \STATE $S^{(0)}_k \leftarrow S_k ,$ for all $k$
   \FOR{$t=1$ {\bfseries to} $\log n$ }
     \STATE $S^{(t)}_k \leftarrow \emptyset ,$ for all $k$
      \FOR{$v \in \mathcal{V}$}
  \STATE Find $k^{\star} = \arg \max_{1\le k\le \hat{K}} \{\sum_{i=1}^{\hat{K}} \sum_{w \in S^{(t-1)}_i}\sum_{\ell =0}^L A_{vw}^\ell \log \hat{p}(k,i,\ell)  \} $ (tie broken uniformly at random)
   \STATE $S^{(t)}_{k^{\star}} \leftarrow S^{(t)}_{k^{\star}} \cup \{ v \}$
   \ENDFOR
   \ENDFOR
\STATE $\hat{\mathcal V}_k \leftarrow S_k^{(\log n)}$, for all $k$
   \STATE {\bfseries Output:} $(\hat{\mathcal V}_k)_{k=1,\ldots,\hat{K}}$.
\end{algorithmic}
\end{algorithm}

\begin{algorithm}[htb]
   \caption{Spectral decomposition }
   \label{alg:specg}
\begin{algorithmic}
 \STATE {\bfseries Input:} $A_{\Gamma}, \tilde{p}$ \\
\STATE {\bf 1. Iterative Power Method with singular value thresholding}
\STATE (Initialization)  $\chi \leftarrow n$, $k \leftarrow 0$, and $\hat{U} \leftarrow 0^{n \times 1}$
\WHILE{$\chi \ge \sqrt{n\tilde{p}}\log (n\tilde{p})$}
\STATE $k\leftarrow k+1$,~~~~~ $U_0 \leftarrow$ $n \times 1$ Gaussian random vector
\STATE (Iterative power method) $U_t \leftarrow (A_{\Gamma})^{\lceil2 \log (n)\rceil} U_{0}$
\STATE (Orthonormalizing $U_t$) $\hat{U}_k \leftarrow \frac{U_t - \hat{U}_{1:k-1} (\hat{U}_{1:k-1}^\top U_t)}{\|U_t - \hat{U}_{1:k-1} (\hat{U}_{1:k-1}^\top U_t)\|_2}$ 
\STATE (Estimating the $k$-th singular value) $\chi \leftarrow \|A_{\Gamma} \hat{U}_k \|_2$
\ENDWHILE
\STATE $\tilde{K} \leftarrow k-1$,~~~~~~ $\hat{V} \leftarrow \hat{U}_{1:\tilde{K}}^\top A_\Gamma $
\STATE {\bf 2. Clustering }
\STATE $\mathcal{V}_R \leftarrow$ a subset of $\Gamma$ obtained by randomly selecting $\lceil \log (n) \rceil$ items of $\Gamma$
\STATE $Q_{v} \leftarrow \{ w \in \Gamma :\| \hat{V}_w  -\hat{V}_v\|_2^2 \le  \frac{n\tilde{p}^2 }{\log(n\tilde{p})} \} $ for all $v\in \mathcal{V}_R$
\STATE (Initialization) $S_{0}\leftarrow \emptyset$, \quad $k\leftarrow 0$, and \quad $\rho \leftarrow |\Gamma|$
\WHILE{$\rho \ge \frac{\log (n\tilde{p})^4}{\tilde{p}}$ }
\STATE $k\leftarrow k+1$,~~~~~ $v_k^{\star} \leftarrow \arg \max_{v\in \mathcal{V}_R} | Q_{v}\setminus \bigcup_{l=0}^{k-1} S_{l} |$,~~~~~ $S_{k} \leftarrow Q_{v_k^{\star}} \setminus \bigcup_{l=0}^{k-1} S_{l} $ and $\rho \leftarrow |S_{k}| .$
\ENDWHILE
\STATE $\hat{K} \leftarrow k-1$
\FOR{$v \in \Gamma\setminus \bigcup_{k=1}^{\hat{K}} S_{k}$}
\STATE $k_{\star} \leftarrow \arg\min_k \| \hat{V}_{v_k^\star}  -\hat{V}_v\|_2$,~~~~~~ $S_{k_{\star}} \leftarrow S_{k_{\star}} \cup \{ v\}$
\ENDFOR
\STATE {\bfseries Output:} $(S_k)_{k=1,\ldots,\hat{K}}$.
\end{algorithmic}
\end{algorithm}

\section{Properties of the divergence $D(\alpha,p)$ and related quantities}

In this section, we prove the two claims of Section 2, as well as other results on the divergence $D(\alpha,p)$ that will be instrumental in the proofs of Theorems.

\subsection{Proof of Claim~\ref{thm:ch}}\label{sec:proof-theor-refthm:c}

$D_{L+}(p(i),p(j))$ is the minimum of the objective function of the following convex optimization problem:
\begin{align}
\min_{y\in \mathcal{P}^{K\times(L+1) }}&\quad\sum_{k=1}^K \alpha_k\left(\sum_{\ell =1}^L   y(k,\ell)\log\left(\frac{y(k,\ell)}{p(i,k,\ell)}\right)+(1-\sum_{\ell=1}^L y(k,\ell))\log\left(\frac{1-\sum_{\ell=1}^L y(k,\ell)}{1-\sum_{\ell=1}^L p(i,k,\ell)}\right)\right)\cr
\mbox{s.t.}&\quad\sum_{k=1}^K \alpha_k KL(y(k), p(i,k)) \ge \sum_{k=1}^K \alpha_k  KL(y(k), p(j,k)).\label{eq:ch}
\end{align}
Note that we define $y(k,0) = 1- \sum_{\ell=1}^L y(k,\ell)$ for all $k$. Since $\bar{p}= o(1)$, one can easily check that the solution of \eqref{eq:ch} has to be $\sum_{\ell=1}^L y(k,\ell) = o(1)$ for all $k$. The objective function converges to infinity when $\sum_{\ell=1}^L y(k,\ell) = \Omega(1)$, while it has $o(\bar{p})$ when $y(k,\ell) = p(j,k,\ell)$ for all $k$ and $\ell$. Thus, we consider $\sum_{\ell=1}^L y(k,\ell) = o(1)$.
The associated Lagrangian is:
\begin{align}
g(y,\lambda)=&\sum_{k=1}^K \alpha_k\left(\sum_{\ell =1}^L   y(k,\ell)\log\left(\frac{y(k,\ell)}{p(i,k,\ell)}\right)+(1-\sum_{\ell=1}^L y(k,\ell))\log\left(\frac{1-\sum_{\ell=1}^L y(k,\ell)}{1-\sum_{\ell=1}^L p(i,k,\ell)}\right)\right)
+  \cr
&\sum_{k=1}^K \alpha_k \lambda \left(\sum_{\ell=1}^{L}
y(k,\ell)\log\left(\frac{p(i,k,\ell)}{p(j,k,\ell)}\right)+(1-\sum_{\ell=1}^{L}y(k,\ell))\log\left(\frac{1-\sum_{\ell=1}^{L}p(i,k,\ell)}{1-\sum_{\ell=1}^{L} p(j,k,\ell)}\right)\right).\label{eq:la}
\end{align}
The derivative of $g(y,\lambda)$ w.r.t. $y(k,\ell)$ is computed as follows:
\begin{align*}
\frac{\partial g(y,\lambda)}{\partial y(k,\ell)} = &\alpha_k \left(
\log\left(\frac{y(k,\ell)}{p(i,k,\ell)}\right)-\log\left(\frac{1-\sum_{m=1}^{L} y(k,m)}{1-\sum_{m=1}^{L} p(i,k,m)}\right)\right)
+  \cr
& \alpha_k \lambda \left(
\log\left(\frac{p(i,k,\ell)}{p(j,k,\ell)}\right)-\log\left(\frac{1-\sum_{m=1}^{L} p(i,k,m)}{1-\sum_{m=1}^{L}p(j,k,m)}\right)\right).
\end{align*}
Observe that, since (A1) holds, $\bar{p}=o(1)$ and $\sum_{\ell=1}^L y(k,\ell) = o(1)$, as $n$ grows large, $\log\left(\frac{1-\sum_{m=1}^L y(k,m)}{1-\sum_{m=1}^L p(i,k,m)}\right)$ and
$\log\left(\frac{1-\sum_{m=1}^L p(i,k,m)}{1-\sum_{m=1}^L p(j,k,m)}\right)$ converges to 0. Thus, \eqref{eq:la} is minimized at
\begin{equation}y(k,\ell) = p(i,k,\ell)\left( \frac{p(j,k,\ell)}{p(i,k,\ell)}\right)^\lambda (1+o(1)).\label{eq:y}\end{equation}
When we put \eqref{eq:y} onto \eqref{eq:la} and use the approximation
$\lim_{x\to 0}\log(1+x) =
x$ (again using $\bar{p}=o(1)$),
\begin{eqnarray*}
&&\min_{y\in \mathcal{P}^{K\times \{0,1 \}}} g(y,\lambda)\cr
&=&\min_{y\in \mathcal{P}^{K\times \{0,1 \}}} \sum_{k=1}^K\sum_{\ell=1}^L \alpha_k \left( o(\bar{p}) + \right.\cr
&&\left. (1-\sum_{\ell=1}^L y(k,\ell))\log\left(\frac{1-\sum_{\ell=1}^L y(k,\ell)}{1-\sum_{\ell=1}^L p(i,k,\ell)}\right)\left(\frac{1-\sum_{\ell=1}^{L}p(i,k,\ell)}{1-\sum_{\ell=1}^{L} p(j,k,\ell)}\right)^\lambda \right) 
\cr
&=&\min_{y\in \mathcal{P}^{K\times \{0,1 \}}} \sum_{k=1}^K\sum_{\ell=1}^L \alpha_k \left( o(\bar{p}) - \right.\cr
&&\left.\sum_{\ell=1}^L y(k,\ell) (1+o(1)) + (1-\lambda)\sum_{\ell=1}^L p(i,k,\ell) (1+o(1)) + \lambda\sum_{\ell=1}^{L} p(j,k,\ell)(1+o(1)) \right) .
\end{eqnarray*}
Therefore, the minimum value of \eqref{eq:ch} is equivalent to
$$\max_{\lambda\in [0,1]}
\sum_{k=1}^K\sum_{\ell=1}^L \alpha_k \left((1-\lambda) p(i,k,\ell)+\lambda p(j,k,\ell) - p(i,k,1)^{1-\lambda}p(j,k,\ell)^\lambda \right) + o(\bar{p}).$$

\subsection{Proof of Claim~\ref{thm:di}}
When $\bar{p} = o(1)$, for all $i\neq j$, $\alpha_i = \frac{1}{K}$ , $p(i,i,\ell) = p(\ell)$, and $p(i,j,\ell)=q(\ell)$, from Claim~\ref{thm:ch},
\begin{eqnarray}
D_{L+} (\alpha,p(i),p(j)) &=&\max_{\lambda \in [0,1]} \sum_{k=1}^K \sum_{\ell=1}^L \alpha_k \left((1-\lambda) p(i,k,\ell)+\lambda p(j,k,\ell) - p(i,k,\ell)^{1-\lambda}p(j,k,\ell)^\lambda \right) \cr
&=&\frac{1}{K} \max_{\lambda \in [0,1]}  \sum_{\ell=1}^L \left(p(\ell)+ q(\ell) - p(\ell)^{1-\lambda}q(\ell)^\lambda - p(\ell)^{\lambda}q(\ell)^{1-\lambda} \right) \cr
&= &\frac{1}{K}  \sum_{\ell=1}^L \left(p(\ell)+ q(\ell) - 2\sqrt{p(\ell)q(\ell)} \right). \label{eq:11}
\end{eqnarray}

Now, since $\sqrt{1+x} = 1+\frac{x}{2}(1+o(1))$ and $\log(1+x) = x(1+o(1))$ when $x=o(1)$,
\begin{eqnarray}
- \frac{2}{K} \log \left(\sum_{\ell=0}^L \sqrt{p(\ell)q(\ell)} \right) &=& - \frac{2}{K} \log \left(\sqrt{p(0)q(0)}+\sum_{\ell=1}^L \sqrt{p(\ell)q(\ell)} \right)\cr
&=& - \frac{2}{K} \log \left(1- \frac{\sum_{\ell=1}^L p(\ell)+q(\ell)}{2}(1+o(1)) +\sum_{\ell=1}^L \sqrt{p(\ell)q(\ell)} \right)\cr
&=&  \frac{2}{K}  \left( \frac{\sum_{\ell=1}^L p(\ell)+q(\ell)}{2} -\sum_{\ell=1}^L \sqrt{p(\ell)q(\ell)} \right) (1+o(1)).\label{eq:12}
\end{eqnarray}

The claim follows from \eqref{eq:11} and \eqref{eq:12}.

\subsection{Other properties}

\begin{lemma}
Let $(i^\star, j^\star) = \arg\min_{i,j} D_{L+}(p(i),p(j))$ and $i^\star < j^\star$. Then, there exists $q\in \mathcal{P}^{K\times (L+1)}$ such that
$$ D(\alpha,p) = \sum_{k=1}^K \alpha_k KL(q(k),p(i^\star,k))  = \sum_{k=1}^K \alpha_k KL(q(k),p(j^\star,k)).$$\label{lem:dq}
\end{lemma}
{\em Proof.}
We check by contradiction that such a $q$ exists. Indeed, assume that 
$$D(\alpha,p) = \sum_{k=1}^K \alpha_k KL(q(k),p(i^\star,k)) > \sum_{k=1}^K \alpha_k KL(q(k),p(j^\star,k)).$$ Then there exists $k_0$ such that $KL(q(k_0),p(i^\star,k_0)) > KL(q(k_0),p(j^\star,k_0))$. Observe that by positivity of the $KL$ divergence, $q(k_0)\neq p(i^\star,k_0)$. Hence by continuity of the $KL$ divergence, we can construct $q'$ such that $q(k)=q'(k)$ for all $k\neq k_0$, and such that: $KL(q(k_0),p(i^\star,k_0))-\epsilon< KL(q'(k_0),p(i^\star,k_0)) < KL(q(k_0),p(i^\star,k_0))$ and $KL(q'(k_0),p(j^\star,k_0)) < KL(q(k_0),p(j^\star,k_0))+\epsilon$ for some $0<\epsilon < (KL(q(k_0),p(i^\star,k_0)) - KL(q(k_0),p(j^\star,k_0)))/2$. With this choice of $q'$, we get:
$$
D(\alpha,p) > \sum_{k=1}^K \alpha_k KL(q'(k),p(i^\star,k)) > \sum_{k=1}^K \alpha_k KL(q'(k),p(j^\star,k)),
$$ 
which contradicts the definition of $D(\alpha,p)$. 
\QED

\begin{lemma}
When $\bar{p}=o(1)$, $$\lim_{n\to \infty} \frac{D(\alpha,p)}{\sum_{k=1}^K \frac{\alpha_k}{2} \left( \sum_{\ell=1}^L (\sqrt{p(i^\star,k,\ell)} - \sqrt{p(j^\star,k,\ell)})^2\right)} \ge 1 .$$ \label{lem:dpq}
\end{lemma}
{\em Proof.} 
Let $(i^\star, j^\star) = \arg\min_{i,j} D_{L+}(\alpha,p(i),p(j))$ and $i^\star < j^\star$. From Lemma~\ref{lem:dq}, there exists $q$ satisfying that 
$$ D(\alpha,p) = \sum_{k=1}^K \alpha_k KL(q(k),p(i^\star,k))  = \sum_{k=1}^K \alpha_k KL(q(k),p(j^\star,k)).$$
Then,
\begin{eqnarray*}
nD(\alpha,p) &=& n\frac{\sum_{k=1}^K \left(\alpha_k KL(q(k),p(i^\star,k))+\alpha_k KL(q(k),p(j^\star,k))\right)}{2} \cr
&=&-n\sum_{k=1}^K \alpha_k \sum_{\ell=0}^L q(k,\ell)\log\left(\frac{\sqrt{p(i^\star,k,\ell)p(j^\star,k,\ell)}}{q(k,\ell)}\right)\cr
&\ge &n\sum_{k=1}^K \alpha_k \sum_{\ell=0}^L \left( q(k,\ell)- \sqrt{p(i^\star,k,\ell)p(j^\star,k,\ell)}\right) \cr
&= &n\sum_{k=1}^K \alpha_k \left( \frac{\sum_{\ell=1}^L (p(i^\star,k,\ell)+p(j^\star,k,\ell))}{2}- \sum_{\ell=1}^L\sqrt{p(i^\star,k,\ell)p(j^\star,k,\ell)}\right)\left(1-o(1)\right)\cr
&=&n\sum_{k=1}^K \frac{\alpha_k}{2} \left( \sum_{\ell=1}^L (\sqrt{p(i^\star,k,\ell)} - \sqrt{p(j^\star,k,\ell)})^2\right)\left(1-o(1)\right).
\end{eqnarray*}
\QED

\begin{lemma}
Under condition (A1), when $\bar{p}=o(1)$,
$\lim\sup_{n \to \infty}\frac{D(\alpha,p)}{\eta \bar{p} L} \le  1.$ \label{lem:dpu}
\end{lemma}
{\em Proof.} From the definition of $D(\alpha,p)$, for any $i\neq j$,
\begin{eqnarray*}
D(\alpha,p) &\le & \max\left\{ \sum_{k=1}^K \alpha_k KL(p(i,k),p(i,k)), \sum_{k=1}^K \alpha_k KL(p(i,k),p(j,k))\right\} \cr
& =& \sum_{k=1}^K \alpha_k KL(p(i,k),p(j,k)) \cr
 &\le & \sum_{k=1}^K \alpha_k \sum_{\ell=1}^L \frac{(p(i,k,\ell)-p(j,k,\ell))^2}{p(j,k,\ell)}(1+o(1))\cr
&\le & \sum_{k=1}^K \alpha_k \sum_{\ell=1}^L \eta \bar{p}(1+o(1))\cr
& = & \eta \bar{p} L(1+o(1)),
\end{eqnarray*}
where we use $\log(1+x) = x(1+o(1))$ when $x= o(1)$.
\QED

\section{Proof of Theorem~\ref{thm:lower}}

The proof consists in an appropriate {\it change-of-measure} argument. The originality of the proof stems from the fact that the change of measures is obtained by a judicious coupling argument \cite{lindvall2002lectures}. In the following, we refer to $\Phi$ as the true stochastic model under which all the observed random labels are generated, and denote by $\mathbb{P}_\Phi=\mathbb{P}$ (resp. $\mathbb{E}_\Phi[\cdot]=\mathbb{E}[\cdot]$) the corresponding probability measure (resp. expectation). We recall that $\Phi$ is defined by the parameters $(\alpha, p)$, and that under $\Phi$, the nodes are first attached to the various clusters according to the distribution $\alpha$, and the labels between two nodes are then generated using distributions $p$. The proof consists in constructing a perturbed stochastic model $\Psi$ coupling the labels generated under $\Phi$ with those generated under $\Psi$. We denote by $\mathbb{P}_\Psi$ (resp. $\mathbb{E}_\Psi[\cdot]=\mathbb{E}[\cdot]$) the probability measure (resp. expectation) under the perturbed model $\Psi$. We then relate the proportion of misclassified nodes under any given clustering algorithm $\pi$ to the distribution under $\mathbb{P}_\Psi$ of a quantity ${\cal Q}$ that resembles the log-likelihood ratio of the observed labels under $\mathbb{P}_\Phi$ and $\mathbb{P}_\Psi$. The analysis of the likelihood ratio finally provides the desired lower bound on the expected misclassified nodes under $\pi$. Next, we detail each step of the proof.   

\medskip
\noindent
{\bf Coupling and the perturbed stochastic model $\Psi$.} Let $(i^\star, j^\star) = \arg\min_{i,j: i<j} D_{L+}(p(i),p(j))$, and let $v^\star$ denote the smallest node index that belongs to cluster $i^\star$ or $j^\star$. If both $\mathcal{V}_{i^\star}$ and $\mathcal{V}_{j^\star}$ are empty, we define $v^\star =n$. Let $q\in [0,1]^{K\times (L+1)}$ satisfy:
$$ 
D(\alpha,p) = \sum_{k=1}^K \alpha_k KL(q(k),p(i^\star,k))  = \sum_{k=1}^K \alpha_k KL(q(k),p(j^\star,k)).
$$
There exists such a $q$ from Lemma~\ref{lem:dq}. Now to define the perturbed stochastic model $\Psi$, we couple the generation of labels under $\Phi$ and $\Psi$ as follows.
\begin{enumerate}
\item We first generate construct the random clusters $\mathcal{V}_1,\ldots, \mathcal{V}_K$ under $\Phi$, and extract $i^\star$, $j^\star$, and $v^\star$. The clusters generated under $\Psi$ are the same as those generated under $\Phi$. For any $v\in \mathcal{V}$, we denote by $\sigma(v)$ the cluster of node $v$. 
\item For all nodes $v,w\neq v^\star$, the labels generated under $\Psi$ are the same as those generated under $\Phi$, i.e., the label $\ell$ is observed on the edge $(v,w)$ with probability $p(\sigma(v),\sigma(w),\ell)$.
\item Under $\Psi$, for any $v\neq v^\star$, the observed label on the edge $(v,v^\star)$ under $\Psi$ is $\ell$ with probability $q(\sigma(v),\ell)$.  
\end{enumerate} 

\medskip
\noindent
{\bf The log-likelihood ratio and its connection to the expected number of misclassified nodes.} Let $x_{v,w}$ denote the label observed on the edge $(v,w)$. We introduce $\mathcal{Q}$, referred to as the pseudo-log-likelihood ratio of the observed labels under $\mathbb{P}_\Phi$ and $\mathbb{P}_{\Psi}$) as: 
\begin{equation}\label{eq:defL}
\mathcal{Q} = \sum_{v=1}^{v^{\star}-1} \log \frac{q(\sigma(v),x_{v^{\star},v})}{p(\sigma(v^{\star}),
\sigma(v),x_{v^{\star},v})}+\sum_{v=v^{\star}+1}^n \log \frac{q(\sigma(v),x_{v^{\star},v})}{p(\sigma(v^{\star}),\sigma(v),x_{v^{\star},v})}.
\end{equation}

Let $\pi$ denote a clustering algorithm with output $(\hat{\mathcal{V}}_k)_{1\le k \le K}$, and let $\set{E} = \bigcup_{1\le k \le K}\hat{\mathcal{V}}_k \setminus \mathcal{V}_k$ be the set of misclassified nodes under $\pi$. Note that in general in our proofs, we always assume without loss of generality that $| \bigcup_{1\le k \le K}\hat{\mathcal{V}}_k \setminus \mathcal{V}_k | \le |\bigcup_{1\le  k \le K} \hat{\mathcal{V}}_{\gamma(k)} \setminus \mathcal{V}_{k}|$ for any permutation $\gamma$, so that the set of misclassified nodes is really $\set{E}$. We denote by $\varepsilon^\pi (n) = |\set{E}|$. Since under $\Phi$, nodes are interchangeable (remember that nodes are assigned to the various clusters in an i.i.d. manner), we have:
$$
n \mathbb{P}_{\Phi}\{v\in \set{E} \}= \mathbb{E}_{\Phi}[\varepsilon^\pi (n)] = \mathbb{E}[\varepsilon^\pi (n)].
$$

Next, we establish a relationship between $\Ex [\varepsilon^\pi (n) ]$ and the distribution of $\mathcal{Q}$ under $\mathbb{P}_\Psi$. For any function $f(n)$, we have:
\begin{equation} \label{eq:8-1}
\mathbb{P}_{\Psi} \{ \mathcal{Q} \le  f(n) \} = \mathbb{P}_{\Psi} \{ \mathcal{Q} \le f(n) , v^{\star}\in \set{E} \} +\mathbb{P}_{\Psi} \{ \mathcal{Q} \le f(n) , v^{\star} \notin \set{E} \}. 
\end{equation}
Using $\mathcal{Q}$, we get:
\begin{eqnarray}\mathbb{P}_{\Psi} \{ \mathcal{Q} \le  f(n) , v^{\star}\in \set{E}\} &= &
  \int_{\{ \mathcal{Q} \le  f(n) , v^{\star}\in \set{E}\}} d\mathbb{P}_{\Psi} \cr
&= & \int_{\{ \mathcal{Q} \le  f(n) , v^{\star}\in \set{E}\}} \exp (\mathcal{Q}) d\mathbb{P}_{\Phi} \cr
&\le & \exp (f(n))\mathbb{P}_{\Phi} \{ \mathcal{Q} \le  f(n) , v^{\star}\in \set{E}\}\cr 
&\le& \exp (f(n))\mathbb{P}_{\Phi} \{ v^{\star}\in \set{E}\}\cr
&\le &\exp (f(n))\frac{\Ex_{\Phi} [ \varepsilon^\pi (n)]}{(\alpha_{i^\star}+ \alpha_{j^\star})n},
\label{eq:bdl-1} \end{eqnarray}
where the last inequality is obtained from the fact that we cannot distinguish between $v^{\star}$ and any other $v \in \mathcal{V}_{\sigma(v^{\star})}$. Indeed,
\begin{eqnarray*}
\mathbb{P}_{\Phi}\{v^\star \in \set{E} \} &=& \mathbb{P}_{\Phi}\{v \in \set{E} | v \in \mathcal{V}_{i^\star} \cup \mathcal{V}_{j^\star}\} \cr
&=& \frac{\mathbb{P}_{\Phi}\{v \in \set{E} , v \in \mathcal{V}_{i^\star} \cup \mathcal{V}_{j^\star}\}}{\mathbb{P}_{\Phi}\{v \in \mathcal{V}_{i^\star} \cup \mathcal{V}_{j^\star}\}}\cr
&\le& \frac{\mathbb{P}_{\Phi}\{v \in \set{E} \}}{\mathbb{P}_{\Phi}\{v \in \mathcal{V}_{i^\star} \cup \mathcal{V}_{j^\star}\}}~=~\frac{\mathbb{E}_{\Phi}[\varepsilon^\pi (n)]}{(\alpha_{i^\star}+\alpha_{j^\star})n}.
\end{eqnarray*}
Furthermore, since under the stochastic model $\Psi$, the observed labels do not depend on whether $v^\star$ belongs to cluster $i^\star$ or $j^\star$, we have:
\begin{eqnarray*} \mathbb{P}_{\Psi} \{  v^{\star} \in \hat{\mathcal{V}}_{i^\star} |v^{\star} \in \mathcal{V}_{i^\star}  \}& =&  \mathbb{P}_{\Psi} \{  v^{\star} \in \hat{\mathcal{V}}_{i^\star} | v^{\star} \in \mathcal{V}_{j^\star}  \}\quad\mbox{and} \cr
\mathbb{P}_{\Psi} \{  v^{\star} \in \hat{\mathcal{V}}_{j^\star} | v^{\star} \in \mathcal{V}_{i^\star}  \}& =&  \mathbb{P}_{\Psi} \{  v^{\star} \in \hat{\mathcal{V}}_{j^\star} | v^{\star} \in \mathcal{V}_{j^\star}  \}.
\end{eqnarray*}
Finally, since $\mathbb{P}_{\Psi} \{  v^{\star} \in \hat{\mathcal{V}}_{i^\star} |v^{\star} \in \mathcal{V}_{i^\star}  \} + \mathbb{P}_{\Psi} \{  v^{\star} \in \hat{\mathcal{V}}_{j^\star} |v^{\star} \in \mathcal{V}_{i^\star}  \} \le 1$, we also have: 
\begin{eqnarray}&&\mathbb{P}_{\Psi} \{ \mathcal{Q} \le  f(n) , v^{\star} \notin \set{E} \}\cr 
&\le&  \mathbb{P}_{\Psi} \{  v^{\star} \notin \set{E} \}  \cr
&=&\frac{\alpha_{i^{\star}}}{\alpha_{i^{\star}}+\alpha_{j^{\star}}} \mathbb{P}_{\Psi} \{  v^{\star} \in \hat{\mathcal{V}}_{i^\star} |v^{\star} \in \mathcal{V}_{i^\star} \} + \frac{\alpha_{j^{\star}}}{\alpha_{i^{\star}}+\alpha_{j^{\star}}} \mathbb{P}_{\Psi} \{  v^{\star} \in \hat{\mathcal{V}}_{j^\star} |v^{\star} \in \mathcal{V}_{j^\star} \}\cr
&=&\frac{\alpha_{i^{\star}}}{\alpha_{i^{\star}}+\alpha_{j^{\star}}} \mathbb{P}_{\Psi} \{  v^{\star} \in \hat{\mathcal{V}}_{i^\star} |v^{\star} \in \mathcal{V}_{i^\star} \} + \frac{\alpha_{j^{\star}}}{\alpha_{i^{\star}}+\alpha_{j^{\star}}} \mathbb{P}_{\Psi} \{  v^{\star} \in \hat{\mathcal{V}}_{j^\star} |v^{\star} \in \mathcal{V}_{i^\star} \}\cr
&\le &\frac{\alpha_{j^\star}}{\alpha_{i^\star}+\alpha_{j^\star}}.
\label{eq:bdlf-1}
\end{eqnarray} 
Combining \eqref{eq:8-1},~\eqref{eq:bdl-1}, and \eqref{eq:bdlf-1}, we conclude that: 
\begin{equation}
\mathbb{P}_{\Psi} \{ \mathcal{Q} \le  f(n) \} \le
\exp(f(n))\frac{\Ex_{\Phi}[\varepsilon^\pi(n) ]}{(\alpha_{i^\star}+\alpha_{j^\star}) n} +
\frac{\alpha_{j^\star}}{\alpha_{i^\star}+\alpha_{j^\star}}.\label{eq:2v1} 
\end{equation}
The previous equation provides the desired generic relationship between $\Ex_{\Phi}[\varepsilon^\pi(n) ]$ and $\mathbb{P}_{\Psi} \{ \mathcal{Q} \le  f(n) \}$ from which can deduce a necessary condition for $\mathbb{E}[\varepsilon^\pi(n)]\le s$. Applying (\ref{eq:2v1}) with \\$f(n) = \log\left(n/\Ex_{\Phi}[\varepsilon^\pi(n) ] \right) - \log(2/\alpha_{i^\star} )$, we have:
\begin{equation}
\mathbb{P}_{\Psi} \{ \mathcal{Q} \le  \log\left(n/\Ex_{\Phi}[\varepsilon^\pi(n) ] \right) - \log(2/\alpha_{i^\star} ) \} \le 1-\frac{\alpha_{i^\star}}{2} < 1- \frac{\alpha_{i^\star}}{4}.\label{eq:2}
\end{equation}
In addition, from Chebyshev's inequality,
\begin{equation}
\mathbb{P}_{\Psi} \left\{ \mathcal{Q} \le  \mathbb{E}_{\Psi}[\mathcal{Q}] + \sqrt{\frac{4}{\alpha_{i^\star}} \mathbb{E}_{\Psi}[(\mathcal{Q}-\mathbb{E}_{\Psi}[\mathcal{Q}])^2]} \right\} \ge 1-\frac{\alpha_{i^\star}}{4}.\label{eq:5}
\end{equation}
From \eqref{eq:2} and \eqref{eq:5}, we deduce that:
$$
\log\left(n/\Ex_{\Phi}[\varepsilon^\pi(n) ] \right) - \log(2/\alpha_{i^\star} ) \le \mathbb{E}_{\Psi}[\mathcal{Q}] + \sqrt{\frac{4}{\alpha_{i^\star}} \mathbb{E}_{\Psi}[(\mathcal{Q}-\mathbb{E}_{\Psi}[\mathcal{Q}])^2]},
$$
and thus, a necessary condition for $\mathbb{E}[\varepsilon^\pi(n)]\le s$ is:
\begin{equation}
\log\left(n/s \right) -
\log(2/\alpha_{i^\star} ) \le \mathbb{E}_{\Psi}[\mathcal{Q}] + \sqrt{\frac{4}{\alpha_{i^\star}} \mathbb{E}_{\Psi}[(\mathcal{Q}-\mathbb{E}_{\Psi}[\mathcal{Q}])^2]}. \label{eq:fnl}
\end{equation}

\medskip
\noindent
{\bf Analysis of the log-likelihood ratio.} In view of (\ref{eq:fnl}), we can obtain a necessary condition for $\mathbb{E}[\varepsilon^\pi (n)]\le s$ if we evaluate $\mathbb{E}_{\Psi}[\mathcal{Q}]$ and $\mathbb{E}_{\Psi}[(\mathcal{Q}-\mathbb{E}_{\Psi}[\mathcal{Q}])^2]$. 

\medskip
\noindent
(i) We first compute $\mathbb{E}_{\Psi}[\mathcal{Q}]$. Note that in view of the definition of $v^\star$, a node whose index is smaller than $v^\star$ cannot be in $\mathcal{V}_{i^\star}$ or $\mathcal{V}_{j^\star}$, whereas a node whose index $v$ is larger than $v^\star$ can be in any cluster (and the cluster of such a $v$ is drawn according to the distribution $\alpha$ independently of other nodes). This slightly complicates the computation of the expectation of the two sums defining $\mathcal{Q}$ in (\ref{eq:defL}). To circumvent this problem, we can observe that $v^\star$ is rather small, i.e., less $\log(n)^2$ with high probability, and that hence, we can approximate $\mathbb{E}_{\Psi}[\mathcal{Q}]$ by $\mathbb{E}_{\Psi}[\sum_{v=v^{\star}+1}^n \log \frac{q(\sigma(v),x_{v^{\star},v})}{p(\sigma(v^{\star}),\sigma(v),x_{v^{\star},v})}]$, which is itself well-approximated by $nD(\alpha,p)$. More formally, since $\mathbb{P}\{v^\star \le m \} = 1-(1-\alpha_{i^\star}-\alpha_{j^\star})^m$,
\begin{equation} 
\mathbb{P}\{v^\star \le \log(n)^2 \} \ge 1-\frac{1}{n^{4}}.\label{eq:8}
\end{equation}
Hence from condition (A1), \eqref{eq:8}, and the definition of $\mathcal{Q}$,
\begin{eqnarray}
\mathbb{E}_{\Psi}[\mathcal{Q}] &= & \mathbb{P}\{v^\star > \log(n)^2 \} \mathbb{E}_{\Psi}[\mathcal{Q} | v^\star > \log(n)^2] + \mathbb{P}\{v^\star \le \log(n)^2 \} \mathbb{E}_{\Psi}[\mathcal{Q} | v^\star \le \log(n)^2]\cr
&\le& \frac{\log\eta}{n^3}   +  \mathbb{E}_{\Psi}[\mathcal{Q} | v^\star \le \log(n)^2] \cr
&\le& \frac{\log\eta}{n^3}+ \mathbb{E}_{\Psi}\left[\sum_{v=1}^{v^{\star}-1} \log \frac{q(\sigma(v),x_{v^{\star},v})}{p(\sigma(v^{\star}),\sigma(v),x_{v^{\star},v})}| v^\star \le \log(n)^2\right] + nD(\alpha,p)\cr
&\le& \frac{\log\eta}{n^3}+ \mathbb{E}_{\Psi}\left[(v^\star-1) \sum_{k\notin \{i^\star,j^\star \}} \frac{\alpha_k KL(q(k),p(\sigma(v^\star,k)))}{1-\alpha_{i^\star}-\alpha_{j^\star}}| v^\star \le \log(n)^2\right] + nD(\alpha,p)\cr
&\le& \left(n+2 \log (n)^2 \log \eta \right)D(\alpha,p) + \frac{\log\eta}{n^3},\label{eq:9}
\end{eqnarray}
where the last inequlaity stems from the fact that $2 KL(q(i),p(\sigma(v^\star,i))) \log \eta \ge KL(q(j),p(\sigma(v^\star,j)))$ for all $i$ and $j$ from condition (A1).

\medskip
\noindent
(ii) To compute $\mathbb{E}_{\Psi}[(\mathcal{Q}-\mathbb{E}_{\Psi}[\mathcal{Q}])^2]$, we evaluate $\mathbb{E}_{\Psi}[(\mathcal{Q}-nD(\alpha,p))^2|\sigma(v^{\star}) = i^{\star} ] $ and $\mathbb{E}_{\Psi}[(\mathcal{Q}-nD(\alpha,p))^2|\sigma(v^{\star}) = j^{\star} ]$. From condition (A1), \eqref{eq:8}, and the definition of $\mathcal{Q}$, 
\begin{eqnarray*}
&&\mathbb{E}_{\Psi}[(\mathcal{Q}-nD(\alpha,p))^2|\sigma(v^{\star}) = i^{\star} ]\cr
&= &\mathbb{P}\{v^\star \le \log (n)^2 \}\mathbb{E}_{\Psi}[(\mathcal{Q}-nD(\alpha,p))^2|\sigma(v^{\star}) = i^{\star}, v^\star \le \log (n)^2 ] \cr
& & + \mathbb{P}\{v^\star > \log (n)^2 \}\mathbb{E}_{\Psi}[(\mathcal{Q}-nD(\alpha,p))^2|\sigma(v^{\star}) = i^{\star}, v^\star > \log (n)^2 ] \cr
&\le &\mathbb{E}_{\Psi}[(\mathcal{Q}-nD(\alpha,p))^2|\sigma(v^{\star}) = i^{\star}, v^\star \le \log (n)^2 ]\cr
&& + \frac{1}{n^4}\mathbb{E}_{\Psi}[(\mathcal{Q}-nD(\alpha,p))^2|\sigma(v^{\star}) = i^{\star}, v^\star > \log (n)^2 ] \cr
&= & O(n\bar{p}). 
\end{eqnarray*}
To derive the above inequality, we have used:
\begin{align*}
\mathbb{E}_{\Psi}&\left[\left(\sum_{v=v^{\star}+1}^n \left(\log \frac{q(\sigma(v),x_{v^{\star},v})}{p(\sigma(v^{\star}),\sigma(v),x_{v^{\star},v})}-D(\alpha,p)\right)\right)^2|\sigma(v^{\star}) = i^{\star} \right]\cr
 &\qquad= \sum_{v=v^{\star}+1}^n\mathbb{E}_{\Psi}\left[\left(\log\frac{q(\sigma(v),x_{v^\star,v})}{p(i^\star,\sigma(v),x_{v^\star,v})}-D(\alpha,p)\right)^2|\sigma(v^{\star}) = i^{\star}\right]\cr
&\qquad=O(n \bar{p})\quad\mbox{and}\cr
\mathbb{E}_{\Psi}&\left[\left(\sum_{v=1}^{v^{\star}-1} \left(\log \frac{q(\sigma(v),x_{v^{\star},v})}{p(\sigma(v^{\star}),\sigma(v),x_{v^{\star},v})}-D(\alpha,p)\right)\right)^2|\sigma(v^{\star}) = i^{\star} \right]\cr
&\qquad=O(v^\star\bar{p} + (v^\star\bar{p})^2 ),
\end{align*}
where we use (A1) and the fact that every label is generated independently. Using the same approach, we can also conclude that
$\mathbb{E}_{\Psi}[(\mathcal{Q}-nD(\alpha,p))^2|\sigma(v^{\star}) = j^{\star} ] = O(n\bar{p})$. In summary, we have:
\begin{equation}
\mathbb{E}_{\Psi}[(\mathcal{Q}-\mathbb{E}_{\Psi}[\mathcal{Q}])^2] = O(n \bar{p}). \label{eq:lvar}
\end{equation}

We are ready to complete the proof of Theorem \ref{thm:lower}. From \eqref{eq:fnl}, \eqref{eq:9}, \eqref{eq:lvar}, and Lemma~\ref{lem:dpq}, when the expected number of misclassified nodes is less than $s$ (i.e., $\Ex[\varepsilon^\pi(n) ] \le s$ ), we must have:
$$\liminf_{n\to \infty} \frac{nD(\alpha,p)}{\log\left(n/s \right)}\ge 1.$$ \QED

\section{Performance of the SP Algorithm -- Proof of Theorem \ref{thm:algorithms}}

{\bf Notations.} We use the standard matrix norm $\| A\|=\sup\limits_{x:\|
  x\|_2=1}\| Ax\|_2$. We denote by $M^\ell$ the expectation of the matrix of $A^\ell$, i.e., $M^\ell_{u,v} = p(i,j,\ell)$ when $u\in \mathcal{V}_i$ and $v\in \mathcal{V}_j$. Let $M = \sum_{\ell=1}^Lw_\ell M^\ell$. 
We define $A_\Gamma$ to denote the adjacency matrix
obtained after trimming (Step 3 in Algorithm \ref{alg:partition}). For any matrix $R\in \mathbb{R}^{n\times n}$, we define the matrix $R_\Gamma$ the square matrix formed by the lines and columns of $R$ whose indexes are in $\Gamma$. Hence, we can define $A^\ell_\Gamma$, $M^\ell_\Gamma$, and $M_{\Gamma}$ where $\Gamma$ is the set of items obtained after the trimming process (Line 3) in the SP algorithm (when taking the expectation to get for example $M_\Gamma$, we condition on $\Gamma$). We introduce the noise matrices $X^\ell_{\Gamma} = A^\ell_{\Gamma} -
M^\ell_{\Gamma}$ and $X_{\Gamma} = \sum_{\ell=1}^L w_\ell X^\ell_{\Gamma}$.
We also denote by $e(v,S,\ell) = \sum_{w \in S} A^{\ell}_{vw}$ the total
number of item pairs with observed label $\ell$ including the item
$v$ and an item from $S$ and $\mu(v,S,\ell) = \frac{e(v,S,\ell)}{|S|}$
the empirical density of label $\ell$. Let $e(v,S) = \sum_{\ell=1}^L
e(v,S,\ell)$ and $\mu(v,S) = [\mu (v,S,\ell)]_{0\le\ell\le L}$. In what follows, $e(v,\mathcal{V})$ is referred to as the {\it degree} of item $v$ (the number of observed labels different than 0 of pairs of items including $v$). 

\medskip
\noindent
{\bf Outline of the proof.} To analyze the performance of the SP algorithm, we first state preliminary lemmas. Lemma~\ref{lem:upe} is concerned with the concentration of the degree of the various items.  
Lemma~\ref{lem:FOspectral} provides an upper bound of the matrix norm of random noise matrix $X^\ell_{\Gamma}$. From these two lemmas, we analyze the performance of the first part of the SP algorithm, and prove Theorem \ref{thm:spec}. To analyze the second part of the SP algorithm consisting of $\log(n)$ improvement iterations, we introduce an appropriate set of items $H$ such that that $\mathcal{V}\setminus H$ is of cardinality less than $s$ with high probability under the condition that $n D(\alpha,p)  - \frac{n\bar{p}}{\log (n\bar{p})^3} \ge  \log(n/s) + \sqrt{\log(n/s)}$. We further bound the rate of improvement of our cluster estimates in each iteration when restricted to the set of items $H$, and deduce that after $\log(n)$ iterations, no item in $H$ is misclassified.  

\subsection{Preliminary lemmas}

\begin{lemma} For every $v\in \mathcal{V}$ and $c \ge 1$, we have
$$\mathbb{P}\{ e(v,\mathcal{V}) \ge 10 c n\bar{p}L\} \le \exp(-10 c n\bar{p} L) .$$
\label{lem:upe}
\end{lemma}
\noindent
{\em Proof.}  From Markov inequality,
\begin{eqnarray*}
\mathbb{P}\{ e(v,\mathcal{V}) \ge 10n\bar{p}L\}& \le & \inf_{\theta > 0}\frac{
\prod_{k=1}^K \mathbb{E}\left[\exp(\theta e(v,\mathcal{V}_k) )\right] }{\exp ( \theta
                                           10c n\bar{p} L)} \cr
&\le & \inf_{\theta > 0}\frac{
\prod_{k=1}^K \big( 1+ \bar{p}L (\exp(\theta ) -1) \big)^{\alpha_k n} }{\exp ( \theta
                                           10 cn\bar{p} L)} \cr
&\le & \inf_{\theta > 0}\frac{
\prod_{k=1}^K \big( \exp( \bar{p}L (\exp( \theta)  - 1 ) ) \big)^{\alpha_k n} }{\exp ( \theta
                                           10c n\bar{p} L)}  \cr
&\le & \exp(-10 cn\bar{p} L),
\end{eqnarray*}
where we derive the last inequality choosing $\theta =2$. \QED

\begin{lemma}[Lemma 8.5 of \cite{coja2010}] When $e(v,\mathcal{V},\ell) \le \Delta$ for all $v \in \Gamma$, 
with high probability, $$\| X^\ell_{\Gamma} \| = O(\sqrt{n\bar{p} + \Delta}).$$
\label{lem:FOspectral}
\end{lemma}

The proof of Lemma \ref{lem:FOspectral} relies on arguments used in
the spectral analysis of random graphs, see \cite{feige2005spectral} and \cite{coja2010}.

\begin{lemma} For all $v \in {\mathcal{V}}_k$ and $D\ge 0$,
\begin{multline*}\mathbb{P}\left\{ \left(\sum_{i=1}^K |\mathcal{V}_i| KL(\mu (v,\mathcal{V}_i),p(k,i)) \ge nD\right)\cap\bigg( e(v,\mathcal{V}) \le 10\eta n\bar{p}L\bigg)\right\} \\ \le  \exp \left(-nD + KL\log(10\eta Ln\bar{p}) +\frac{100
   \eta^2 n\bar{p}^2L^2}{\alpha_1} \right).\end{multline*}\label{lem:kl1}
\end{lemma}

\noindent
{\em Proof.} Let $\mathcal{X}$ be a set of $K\times (L+1)$ matrices such that
$$\mathcal{X} = \left\{\boldsymbol{x}\in \mathbb{Z}^{K\times (L+1)}:\quad \sum_{i=1}^{K}\sum_{\ell=1}^{L} x_{i,\ell} \le 10 \eta n\bar{p}L ,\quad\mbox{and}\quad \sum_{\ell=0}^{L} x_{i,\ell} = |\mathcal{V}_i| \quad\mbox{for all}\quad 1\le i\le K \right\}. $$
For notational simplicity, we use $[\frac{x_{i,\ell}}{|\mathcal{V}_i|}]$ instead of
$[\frac{x_{i,\ell}}{|\mathcal{V}_i|}]_{0\le \ell \le L}$ to represent the
probability mass vector on labels defined by $x_{i}$. With a slight abuse of notation, we denote by $e(v)$ the $K\times
(L+1)$ matrix whose $(i,\ell)$ element is $e(v,\mathcal{V}_i,\ell)$. Then, for $v \in
\mathcal{V}_k$,
\begin{align*}
\mathbb{P}&\left\{  \left( \sum_{i=1}^K |\mathcal{V}_i| KL(\mu(v,\mathcal{V}_i),p(k,i) ) \ge nD \right)\cap\bigg( e(v,\mathcal{V}) \le 10 n\bar{p}L\bigg) \right\} \cr
= & \ \sum_{\boldsymbol{x}\in \mathcal{X}}
     \mathbb{P}\left\{e(v)=\boldsymbol{x} \right\} \mathbb{P}\left\{ \sum_{i=1}^K |\mathcal{V}_i| KL(\mu(v,\mathcal{V}_i),p(k,i) ) \ge nD \bigg| e(v)=\boldsymbol{x}
     \right\} \cr
\le& \ \sum_{\boldsymbol{x}\in \mathcal{X}}
     \mathbb{P}\{e(v)=\boldsymbol{x} \}\frac{\exp\left(\sum_{i=1}^K |\mathcal{V}_i| KL([\frac{x_{i,\ell}}{|\mathcal{V}_i|}],p(k,i) )\right)}{\exp(n D)} \cr
\le& \ \sum_{\boldsymbol{x}\in \mathcal{X}}
\mathbb{P}\{e(v)=\boldsymbol{x} \}\frac{\prod_{i=1}^{K}\prod_{\ell=0}^{L}\left(\frac{x_{i,\ell}}{|\mathcal{V}_i| p(k,i,\ell)} \right)^{x_{i,\ell}}}{\exp(n D)} \cr
\stackrel{(a)}{\le}
&\ \frac{1}{\exp(n D)} \sum_{\boldsymbol{x}\in \mathcal{X}}
 \prod_{i=1}^K\left( \left(1-\frac {\sum_{\ell=1}^L x_{i,\ell}}{|\mathcal{V}_i| 
   }\right)^{x_{i,0}}\exp( \sum_{\ell =1 }^L x_{i,\ell} ) \right)\cr
=&\ \frac{1}{\exp(n D)} \sum_{\boldsymbol{x}\in \mathcal{X}}
 \prod_{i=1}^K\exp\left((|\mathcal{V}_i|-\sum_{\ell=1}^L x_{i,\ell}) \log\left(1-\frac {\sum_{\ell=1}^L x_{i,\ell}}{|\mathcal{V}_i|
   }\right) +  \sum_{\ell =1 }^L x_{i,\ell}  \right)\cr
\le \ & \frac{1}{{\exp(n D)}}\sum_{\boldsymbol{x}\in \mathcal{X}}
     \prod_{i=1}^K\exp\left( \frac{(\sum_{\ell=1}^L x_{k,\ell})^2}{|\mathcal{V}_i|}
     \right) \cr
\le\ &  \frac{(10\eta n\bar{p}L)^{KL} \exp(100 \eta^2 n\bar{p}^2L^2/\alpha_1)}{\exp(n D)} \cr
=&\  \exp \left(-nD + KL\log(10\eta Ln\bar{p}) +\frac{100
   \eta^2 n\bar{p}^2L^2}{\alpha_1} \right),
\end{align*}
where  $(a)$ stems from the following inequality:
\begin{align*}
\mathbb{P}\{e(v,\mathcal{V}_i,\ell)=x_{i,\ell}&\quad\mbox{for all}\quad
  i,\ell \} \cr
\le&\prod_{i=1}^K\left(  p(k,i,0)^{x_{i,0}}
     \prod_{\ell=1}^L {|\mathcal{V}_i| \choose x_{i,\ell}}
     p(k,i,\ell)^{x_{k,\ell}}\right)\cr
\le&\prod_{i=1}^K\left( p(k,i,0)^{x_{i,0}}
     \prod_{\ell=1}^L \left(\frac{e |\mathcal{V}_i| }{x_{i,\ell}}\right)^{x_{i,\ell}}
     p(k,i,\ell)^{x_{i,\ell}}\right).
\end{align*}
\QED

\subsection{Part 1 of the SP algorithm -- Proof of Theorem~\ref{thm:spec}}\label{sec:part-1-sp}

Recall that $\hat{A} = \hat{U}\hat{V} = \hat{U}\hat{U}^\top A_\Gamma$ and $\|\hat{A}_u -\hat{A}_v \| = \|\hat{V}_u-\hat{V}_v \|$. We can bound the number of misclassified items as follows:
\begin{itemize}
\item with high probability, we have
\begin{equation}
\|\hat{A}-M_\Gamma \|_F^2 =\sum_{v\in \Gamma} \| \hat{A}_v - M_{v,\Gamma} \|_2^2 = O(n\bar{p}\log(n\bar{p})^2); \label{eq:am-noise}
\end{equation}
\item with high probability, every item pair $u$ and $v$ satisfies that when $\sigma(v)$ represents the cluster of $v$ and $M_{v,\Gamma}$ denotes
the column vector of $M_\Gamma$ on $v$,
\begin{equation}\|M_{u,\Gamma}- M_{v,\Gamma} \|_2^2 = \Omega \left( n\bar{p}^2 \right)\quad\mbox{when}\quad \sigma(u)\neq \sigma(v),\label{eq:4}
\end{equation}
since every $w_\ell$ is generated uniformly at random in $[0,1]$ and (A2) holds;
\item \eqref{eq:4} suggests that if $v$ is misclassified by Algorithm~\ref{alg:specg}, then we should have:
\begin{equation}
\| \hat{A}_v - M_{v,\Gamma} \|_2^2 = \Omega \left( n\bar{p}^2 \right);\label{eq:10}
\end{equation}
\item from \eqref{eq:am-noise} and \eqref{eq:10}, with high probability,
$$
\left|\bigcup_{k=1}^{K}(V_{k}\setminus S_k )\right|  = O\left(\frac{\log(n\bar{p})^2}{\bar{p}}\right).
$$
\end{itemize}
Next, we prove \eqref{eq:am-noise} and \eqref{eq:10}. 


\smallskip
\noindent{\em Proof of \eqref{eq:am-noise}.}
First observe that from the definition of $\Gamma$, 
\begin{eqnarray*}
\mathbb{P}\left\{ \max_{v\in \Gamma} e(v,\mathcal{V}) \ge 10n \bar{p} L \right\} & = & \mathbb{P}\left\{ |\{v: e(v,\mathcal{V}) \ge 10n \bar{p} L\}| >  \lfloor n \exp(-n\tilde{p}) \rfloor \right\}\cr
&\le&\frac{n\exp(-10n\bar{p}L)}{\lfloor n \exp(-n\tilde{p}) \rfloor +1 }\cr
&\le&\exp(-5n\bar{p}L),
\end{eqnarray*}
where the first inequality stems from Lemma~\ref{lem:upe} and Markov inequality. Therefore, with high probability,
\begin{equation}
\max_{v\in \Gamma} e(v,\mathcal{V}) \le 10n \bar{p} L. \label{eq:ebnd}
\end{equation}
When the degrees of items are bounded, the standard matrix norm of each noise matrix $X_{\Gamma}^\ell$ can be bounded using Lemma~\ref{lem:FOspectral}.
 From \eqref{eq:ebnd} and Lemma~\ref{lem:FOspectral},
\begin{eqnarray}
\|X_\Gamma \| & \le & \sum_{\ell=1}^L w_\ell \|X_\Gamma^\ell \| \cr
&= & \sum_{\ell=1}^L O(w_\ell \sqrt{n\bar{p} +10n\bar{p}L})\cr
&=& O( \sqrt{n\bar{p}}).\label{eq:part1-noise}
\end{eqnarray}

Let $\tilde{K}$ be the number of columns of $\hat{U}$. Since $\hat{A}$ is the $\tilde{K}$-rank approximation of $A_\Gamma$ obtained by the iterative power method with $2 \log (n)$ iterations,
from Theorem 9.1 and Theorem 9.2 in \cite{halko2011finding}, with high probability,
\begin{equation}
\frac{1}{2} s_{k} (A_{\Gamma}) \le \|A_{\Gamma}\hat{U}_{k} \| \le s_{k} (A_{\Gamma}) \quad\mbox{and}\quad\|A_{\Gamma}(I-\hat{U}_{1:k}\hat{U}_{1:k}^\top ) \| \le 2 s_{k+1} (A_{\Gamma}). \label{eq:part1-ipm}
\end{equation}
Since $\|A_{\Gamma}\hat{U}_{K} \| \le s_{K+1} (A_{\Gamma}) \le \|X_{\Gamma} \| = O(\sqrt{n\bar{p}})$ from Lemma~\ref{lem:FOspectral} and \eqref{eq:part1-ipm}, $\tilde{K} \le K$ and thus the rank of $(\hat{A} - M_\Gamma)$ is less than $2K$. Therefore,
\begin{eqnarray}
\|\hat{A} - M_\Gamma \|_F^2 &\le& 2K \|\hat{A} - M_\Gamma \|^2\cr
&\le& 4K\left( \|\hat{A} - A_\Gamma \|^2 + \|A_\Gamma - M_\Gamma \|^2 \right)\cr
&\le& O(n\bar{p} \log(n\bar{p})^2),\label{eq:3}
\end{eqnarray}
where the last inequality stems from the fact that $\|A_\Gamma - M_\Gamma \| = \|X_\Gamma \| = O( \sqrt{n\bar{p}})$ and $\|\hat{A} - A_\Gamma \| \le 2 s_{\tilde{K}+1} (A_{\Gamma}) = O(\sqrt{n\bar{p}}\log (n\bar{p}))$ from \eqref{eq:part1-ipm}.

\smallskip
\noindent{\em Proof of \eqref{eq:10}.}
Define the following sets: 
\begin{eqnarray*}
\mathcal{I}_k & = &\{ v \in \mathcal{V}_k \cap \Gamma : \|\hat{A}_v - M_{\Gamma}^{k} \|^2 \le \frac{1}{4} \frac{n\tilde{p}^2}{\log(n\tilde{p})} \} \cr
\mathcal{O} & = &\{ v \in \Gamma : \|\hat{A}_v - M_{\Gamma}^{k}
\|^2 \ge 4 \frac{n\tilde{p}^2}{\log(n\tilde{p})}\quad\mbox{for all}\quad 1\le
k\le K \}.
\end{eqnarray*}
These sets are designed so that
\begin{itemize}
\item[(i)] $|(\cup_{k=1}^{K} \mathcal{I}_k) \cap
Q_{v} | =0$ for all $v\in \mathcal{O}\cap \mathcal{V}_R$, since $\|\hat{A}_v - \hat{A}_w \|^2 \ge
\frac{1}{2}\|\hat{A}_v -M_{\Gamma}^{k}\|^2 - \|\hat{A}_w
-M_{\Gamma}^{k}\|^2 > \frac{n\tilde{p}^2}{\log(n\tilde{p})} $ for all $w \in \mathcal{I}_k$;
\item[(ii)] 
$|\Gamma \setminus (\cup_{k=1}^K\mathcal{I}_k) | \le \frac{\|\hat{A} - M_\Gamma \|_F^2}{\min_{v \in \Gamma \setminus (\cup_{k=1}^KI_k)  } \|\hat{A}_v - M_{\Gamma}^{k} \|^2} = O\left(\frac{\log(n\bar{p})^3}{\bar{p}} \right)$; 
\item[(iii)]  $\mathcal{I}_k \subset Q_{v}$ for all $v\in \mathcal{I}_k \cap \mathcal{V}_R$, since $\|\hat{A}_v - \hat{A}_w \|^2 \le
2\|\hat{A}_v -M_{\Gamma}^{k}\|^2 + 2\|\hat{A}_w -M_{\Gamma}^{k}\|^2 \le \frac{n\tilde{p}^2}{\log(n\tilde{p})}$ for all $w \in \mathcal{I}_k$;
\item[(iv)] If $|Q_{v}\cap \mathcal{I}_k | \ge 1 $, $|Q_{v}\cap \mathcal{I}_j | =0 $ for all $j\neq k$, since $\|M_{\Gamma}^k - M_{\Gamma}^j \| = \Omega(n\bar{p}^2)$ is much larger than the radius $ \frac{n\tilde{p}^2}{\log(n\tilde{p})} = O(\frac{n\bar{p}^2}{\log(n\bar{p})})$;
\end{itemize}
From the properties of $\mathcal{I}_k$ and $\mathcal{O}$, we state the following results.
\begin{itemize}
\item {From (i)} and {(ii)}, we deduce that
\begin{equation}
|Q_{v}| = O\left(\frac{\log(n\bar{p})^3}{\bar{p}} \right)\quad \mbox{for all} \quad v \in \mathcal{O}\cap V_R,\label{eq:15}
\end{equation} 
since every $w\in (\cup_{k=1}^K\mathcal{I}_k)$ is outside of $Q_{v}$ (i.e., $w\in \Gamma \setminus (\cup_{k=1}^KI_k)$ is necessary for $w\in Q_{v}$);  

\item since $\alpha_k$ is a constant for all $k$ and $\frac{|\Gamma \setminus (\cup_{k=1}^K\mathcal{I}_k) |}{|\Gamma |} = o(1)$ from {(ii)}, with high probability, 
\begin{equation}
|\mathcal{I}_k \cap \mathcal{V}_R| \ge 1 \quad\mbox{for all}\quad 1\le k\le K ;\label{eq:14}
\end{equation}
\item The properties {(ii)}, {(iii)}, and {(iv)} and \eqref{eq:14} imply that 
\begin{equation}
|Q_{v} \setminus \cup_{l=0}^{k-1} S_{l}| \ge m_k, \quad \exists v \in (\cup_{m=1}^K\mathcal{I}_k\cap V_R ) \setminus (\cup_{l=0}^{k-1} S_{l}),\label{eq:16}
\end{equation}
where $m_k$ is the $k$-th largest value among $\{|\mathcal{I}_1|,\dots,|\mathcal{I}_K|\}$ ;
\item since $ |\mathcal{I}_k| \ge |V_k\cap(\Gamma\setminus \mathcal{O}) | \ge \alpha_k n (1-o(1))$ from (ii) and {(iii)},
\begin{equation}
|\mathcal{I}_k| \ge |V_k\cap(\Gamma\setminus \mathcal{O}) | \ge \alpha_k n (1-o(1)).\label{eq:17}
\end{equation}
\end{itemize}
Thus, we can conclude that $\hat{K} = K$ from \eqref{eq:16} and \eqref{eq:17} and the property (ii); and  from \eqref{eq:15}, there exists a permutation $\gamma$ such that $\|\hat{A}_{v^\star_{k}} -M_{\Gamma}^{\gamma(k)} \|^2 \le 4\frac{n\tilde{p}^2}{\log (n\tilde{p})}$ for all $k$. 
Hence from \eqref{eq:4}, 
$\| \hat{A}_v - M_{v,\Gamma} \|^2 = \Omega \left( n\bar{p}^2 \right)$ when $v$ is misclassified.
\QED


\subsection{Proof of Theorem~\ref{thm:algorithms}}
From Chernoff bound, with high probability, 
\begin{equation}
\left| |\mathcal{V}_k|-\alpha_k n \right| \le  \sqrt{n}\log(n)\quad\mbox{for all}\quad k.\label{eq:size}\end{equation} 
In what follows, we hence just prove the theorem assuming that \eqref{eq:size} holds.

Let $H$ be the largest set of items $v\in \mathcal{V}$ satisfying:
\begin{itemize}
\item[(H1)] $e(v,\mathcal{V}) \le 10 \eta n\bar{p}L$,
\item[(H2)] When $v\in \mathcal{V}_k$, $ \sum_{i =1}^K \sum_{\ell =0}^L e (v,\mathcal{V}_i ,\ell ) \log
  \frac{p (k,i,\ell)}{p (j,i,\ell)}  \ge  \frac{n\bar{p}}{\log (n\bar{p})^4}$ for all $j \neq k$.
\item[(H3)] $e(v,\mathcal{V}\setminus H) \le 2 \log (n\bar{p})^2.$
\end{itemize}
(H1) regularizes degrees, (H2) means that $v\in H$ is correctly classified when using the log-likelihood estimate, and (H3) means that $v$ does not share too many labels with items outside $H$.

The proof of the theorem follows from the following propositions. The first provides an upper bound of $|\mathcal{V}\setminus H|$, and the second provides the rate at which our estimated clusters improve in each iteration when we restrict our attention to items in $H$. 

\begin{proposition}
When $n D(\alpha,p)  - \frac{n\bar{p}}{\log (n\bar{p})^3} \ge  \log(n/s) + \sqrt{\log(n/s)}$, $|\mathcal{V}\setminus H| \le s$ with high probability.\label{prop:1}
\end{proposition}

\begin{proposition} If ${| \bigcup_{k=1}^K (S^{(0)}_k \setminus \mathcal{V}_k)\cap H| + |\mathcal{V} \setminus
  H|} = O(1/\bar{p})$, with high probability, the following statement holds 
$$ \frac{ | \bigcup_{k=1}^K (S^{(t+1)}_k \setminus \mathcal{V}_k)\cap H|}{
  |\bigcup_{k=1}^K (S^{(t)}_k \setminus \mathcal{V}_k)\cap H|} \le \frac{1}{\sqrt{n\bar{p}}}\quad\mbox{for all}\quad t\ge 0.$$ \label{lem:improve}
\end{proposition}

From Proposition \ref{lem:improve}, after $\log(n)$ iterations (remember that $n\bar{p}=\omega(1)$, so when $n$ is large enough $1/\sqrt{n\bar{p}}\le e^{-2}$), no item in $H$ can be misclassified with high probability. Hence the number of misclassified items cannot exceed $|\mathcal{V}\setminus H|\le s$, $n D(\alpha,p)  - \frac{n\bar{p}}{\log (n\bar{p})^3} \ge  \log(n/s) + \sqrt{\log(n/s)}$. The proof is completed by remarking that if the previous condition on $D(\alpha, p)$ holds, then
$$1 \le \lim_{n\to \infty} \frac{n D(\alpha,p)  - \frac{n\bar{p}}{\log (n\bar{p})^3} }{\log(n/s) + \sqrt{\log(n/s)}} = \lim_{n\to \infty} \frac{n D(\alpha,p)}{\log(n/s) },$$
where we used $D(\alpha, p)= \Omega(\bar{p})$ from condition (A2) and Lemma~\ref{lem:dpq}.
\QED

\subsubsection{Proof of Proposition \ref{prop:1} -- Size of $\mathcal{V}\setminus H$}

We compute the number of items satisfying (H1), (H2), and (H3) in \eqref{eq:h1}, \eqref{eq:h2}, and Lemma \ref{lem:sizeH}, respectively. 

\medskip
\noindent
\underline{Number of items satisfying (H1):} From Lemma~\ref{lem:upe}, we get:
\begin{equation}
\mathbb{P}\{ e(v,\mathcal{V}) \le 10\eta n\bar{p}L\} \ge 1- \exp(-10 \eta n\bar{p} L).\label{eq:h1}
\end{equation}

\medskip
\noindent
\underline{Number of items satisfying (H2):} We shall prove that when $v$ satisfies (H1), $v$ satisfies (H2) as well with probability at least
\begin{equation}
1-\exp\left( -nD(\alpha,p)+ \frac{n \bar{p}}{ 2\log (n
\bar{p})^3} \right).\label{eq:h2}
\end{equation}

To this aim, we first establish that if $v$ satisfies
\begin{equation}\sum_{i=1}^K |\mathcal{V}_i| KL(\mu(v,\mathcal{V}_i),p(k,i))  \le \left( 1- \frac{\log (n)^2}{ \sqrt{n}}\right) n D(\alpha,p) - \frac{n \bar{p}}{ \log (n\bar{p})^4},\label{eq:KLD}\end{equation}  
then $v$ satisfies (H2). Indeed, assume that \eqref{eq:KLD} holds, then 
\begin{itemize} 
\item[(i)] $\sum_{i=1}^K \alpha_i n KL(\mu(v,\mathcal{V}_i),p(k,i)) \le \left( 1+ \frac{\log (n)^2}{ \sqrt{n}}\right)\sum_{i=1}^K |\mathcal{V}_i| KL(\mu(v,\mathcal{V}_i),p(k,i)) < n D(\alpha,p)$, since $||\mathcal{V}_i|-\alpha_i n| \le \sqrt{n}\log(n)$ and \eqref{eq:KLD} holds; 
\item[(ii)] $\sum_{i=1}^K \alpha_i n KL(\mu(v,\mathcal{V}_i),p(j,i)) \ge n D(\alpha,p)$, since\\ $\max\left\{\sum_{i=1}^K \alpha_i KL(\mu(v,\mathcal{V}_i),p(j,i)),\sum_{i=1}^K \alpha_i KL(\mu(v,\mathcal{V}_i),p(k,i))\right\} \ge D(\alpha,p)$ and \\$\sum_{i=1}^K \alpha_i KL(\mu(v,\mathcal{V}_i),p(k,i)) < D(\alpha,p)$; 
\item[(iii)] $\sum_{i=1}^K |\mathcal{V}_i| KL(\mu(v,\mathcal{V}_i),p(j,i)) \ge \left( 1 - \frac{\log (n)^2}{ \sqrt{n}}\right)n D(\alpha,p)$, from ii) and the fact that $||\mathcal{V}_i|-\alpha_i n| \le \sqrt{n}\log(n)$;
\item[(iv)] from \eqref{eq:KLD} and iii), for all $j\neq i$,
\begin{eqnarray*}\sum_{i =1}^K \sum_{\ell =0}^L e (v,
  \mathcal{V}_i ,\ell ) \log \frac{p (k,i,\ell)}{p (j,i,\ell)} & =&  \sum_{i=1}^K |\mathcal{V}_i| \left( KL(\mu(v,\mathcal{V}_i),p(j,i)) -KL(\mu(v,\mathcal{V}_i),p(k,i))\right) \cr
&\ge& \frac{n \bar{p}}{ \log (n
\bar{p})^4}.
\end{eqnarray*}
\end{itemize}
Hence $v$ satisfies (H2). It remains to evaluate the probability of the event \eqref{eq:KLD}, which is done by applying Lemma \ref{lem:kl1} and proves \eqref{eq:h2}.

\medskip
\noindent
\underline{Number of items satisfying (H3):} From \eqref{eq:h1}, \eqref{eq:h2}, and the Markov inequality, we deduce that with probability at least $1- \exp \left(-\sqrt{\log(n/s)} \right)$, the number of items that do not satisfy either (H1) or (H2)  is less than $s/3$ when $n D(\alpha,p)  - \frac{n\bar{p}}{\log (n\bar{p})^3} \ge  \log(n/s) + \sqrt{\log(n/s)}$, since
\begin{align}&\frac{\mathbb{E}\{\mbox{The number of items that
      do not satisfy either (H1) or (H2)}\}}{s/3} \cr
&\quad\le \frac{n \exp(-10 \eta n\bar{p} L) + n \exp\left( -nD(\alpha,p)+ \frac{n \bar{p}}{ 2\log (n
\bar{p})^3} \right)}{s/3} \cr
&\quad\le \frac{n}{s} \exp\left( -nD(\alpha,p)+ \frac{n \bar{p}}{ \log (n
\bar{p})^3} \right) \le \exp \left(-\sqrt{\log(n/s)} \right) ,\label{eq:H1H2}\end{align}
where we have used Lemma~\ref{lem:dpu} for the last inequality. Lemma~\ref{lem:sizeH} allows us to complete the proof of Proposition. \QED 

\begin{lemma}
When the number of items that do not satisfy either (H1) or (H2)  is less than $s/3$,
$|\mathcal{V} \setminus H| \le s$, with high probability. 
\label{lem:sizeH}
\end{lemma}

\medskip\noindent
{\em Proof.} Let $e(S,S) = \sum_{v\in S} e(S,S)$. Next we prove the following intermediate claim: there is no subset
$S\subset \mathcal{V}$ such that $e(S,S) \ge s \log (n\bar{p})^2$ and $|S|= s$ with high probability. For any subset $S \in \mathcal{V}$
such that $|S| = s,$ by Markov inequality,
\begin{eqnarray} 
\mathbb{P} \{ e(S,S) \ge s \log (n\bar{p})^2 \} &\le&\inf_{t\ge  0} \frac{ \mathbb{E}[ \exp (e(S,S)t) ]  }{ st \log ( n\bar{p})^2 } \cr
&\le&  \inf_{t\ge 0} \frac{ \prod_{i=1}^{s^2/2} ( 1+ L\bar{p} \exp(t) )  }{ st \log (n\bar{p})^2 } \cr
&\le & \inf_{t\ge 0} \exp \left( \frac{s^2L\bar{p}}{2} \exp (t) - st \log (n\bar{p})^2  \right) \cr
&\le &  \exp \left(   -   n\bar{p}s\big(\log n\bar{p} -
       \frac{sL}{2n}\exp(\frac{n\bar{p}}{\log n\bar{p}})\big) \right)\cr 
&\le& \exp \left(   - \frac{n\bar{p}s \log n\bar{p}}{2} \right),\label{eq:bndss}
\end{eqnarray}
where, in the last two inequalities, we have set $t= \frac{n\bar{p}}{\log n\bar{p}}$ and used the fact that: 
$\frac{n}{s} \ge \exp(\frac{n\bar{p}}{\log n\bar{p}}),$
which comes from the assumptions made in the theorem. 
Since the number of subsets $S \subset \mathcal{V}$ with size $s$ is
${{n}\choose{s}} \le (\frac{e n}{s})^{s} ,$ from \eqref{eq:bndss}, we deduce:
\begin{align*} \mathbb{E}[| \{S : e(S,S) \ge s \log (n\bar{p})^2 ~\mbox{and}~|S|= s \} |]  
&\le (\frac{e n}{s})^{s} \exp \left(   -   \frac{n\bar{p}s \log n\bar{p}}{2}
  \right) \cr
&= \exp \left( -s(\frac{n\bar{p} \log n\bar{p}}{2} - \log \frac{en}{s}) \right) \cr
&\le \exp\left(-\frac{n\bar{p}s \log n\bar{p}}{4}\right).\end{align*}
Therefore, by Markov inequality, we can conclude that there is no $S \subset \mathcal{V}$ such that
$e(S,S) \ge s \log (n\bar{p})^2$ and $|S|= s$ with high probability.

To conclude the proof of the lemma, we build the following sequence of sets. Let $Z_1$ denote the set of items that do not satisfy at least one
 of (H1) and (H2). Let $\{ Z(t) \subset \mathcal{V}\}_{1\le t \le t^{\star}}$ be generated as follows:
\begin{itemize}
\item $Z(0)=Z_1$.
\item For $t \ge 1$, $Z(t) = Z(t-1) \cup \{v_t \}$ if there exists
  $v_t \in \mathcal{V}$ such that
  $e(v_t , Z(t-1)) > 2 \log (n\bar{p})^2$ and
  $v_t \notin Z(t-1)$. If such an item does not exist, the sequence ends.
 \end{itemize}
The sequence ends after the construction of $Z(t^\star)$. We show that if we assume that the cardinality of items that do not satisfy (H3) is strictly larger than $s/2$, then one the set of the sequence $\{ Z(t) \subset \mathcal{V}\}_{1\le t \le t^{\star}}$ contradicts the claim we just proved.

Assume that the number of items do not satisfy (H3) is strictly larger than $s/2$, then these items will be at some point added to the sets $Z(t)$, and by definition, each of these node contributes with more than $2 \log (n\bar{p})^2$ in $e(Z(t),Z(t))$. Hence if starting from $Z_1$, we add $s/2$ items not satisfying (H3), we get a set $Z(t)$ of cardinality less than $s/3+s/2$ and such that $e(Z(t),Z(t))> s\log (n\bar{p})^2$. We can further add arbitrary items to $Z(t)$ so that it becomes of cardinality $s$, and the obtained set contradicts the claim. \QED

\subsubsection{Proof of Proposition~\ref{lem:improve}}
Recall that $\{ S^{(t)}_j\}_{1\le j\le K}$ is the partition after the
$t$-th improvement iteration. Also recall that with loss of generality, we assume that the set of misclassified items in $H$ after the $t$-th step is $\mathcal{E}^{(t)}=\left(\cup_k (S_k^{(t)}\setminus \mathcal{V}_k)\right)\cap H$ (it should be defined through an appropriate permutation $\gamma$ of $\{1,\ldots,K\}$ by $\mathcal{E}^{(t)}=(\cup_k (S_k^{(t)}\setminus \mathcal{V}_{\gamma(k)}))\cap H$, but we omit $\gamma$). With this notational convention, we can define $\mathcal{E}_{jk}^{(t)} = (S^{(t)}_j \cap
\mathcal{V}_k)\cap H$ and $\mathcal{E}^{(t)}= \bigcup_{j,k:j \neq k}
\mathcal{E}_{jk}^{(t)}$.
At each improvement step, items move to the most likely cluster (according to the log-likelihood defined in the SP algorithm). Thus, for all $i$,
\begin{align}
0\le &\sum_{j,k:j\neq k}\sum_{v \in \mathcal{E}_{jk}^{(t+1)}}
       \sum_{i=1}^K\sum_{\ell =0}^L e(v, S^{(t)}_i,\ell) \log
       \frac{\hat{p}(j,i,\ell)}{\hat{p}(k,i,\ell)} \cr
\le &\sum_{j,k:j\neq k}\sum_{v \in \mathcal{E}_{jk}^{(t+1)}}
       \sum_{i=1}^K\sum_{\ell =0}^L e(v, S^{(t)}_i,\ell) \log
       \frac{p(j,i,\ell)}{p(k,i,\ell)} +  |\mathcal{E}^{(t+1)}|(n\bar{p})^{1-\kappa}\log(n\bar{p})^3 \label{eq:7-1} \\
\le &\sum_{j,k:j\neq k}\sum_{v \in \mathcal{E}_{jk}^{(t+1)}}
       \sum_{i=1}^K\sum_{\ell =0}^L e(v, \mathcal{V}_i,\ell) \log
       \frac{p(j,i,\ell)}{p(k,i,\ell)}  \cr
& + \sum_{w\in \mathcal{E}^{(t+1)}} e(w, \mathcal{E}^{(t)})\log (2\eta) + 2|\mathcal{E}^{(t+1)}|(n\bar{p})^{1-\kappa}\log(n\bar{p})^3\label{eq:7-2}\\
\le &-\frac{n\bar{p}}{\log (n\bar{p})^4} |\mathcal{E}^{(t+1)}| +
 \sum_{w\in \mathcal{E}^{(t+1)}} e(w, \mathcal{E}^{(t)},\ell)\log (2\eta) +2|\mathcal{E}^{(t+1)}|(n\bar{p})^{1-\kappa}\log(n\bar{p})^3\label{eq:7-3} \\
\le &-\frac{n\bar{p}}{\log (n\bar{p})^4}|\mathcal{E}^{(t+1)}| + \sqrt{
      |\mathcal{E}^{(t)}||\mathcal{E}^{(t+1)}|n\bar{p} \log n\bar{p}} + 3|\mathcal{E}^{(t+1)}|(n\bar{p})^{1-\kappa}\log(n\bar{p})^3. \label{eq:7-4}
\end{align}
Therefore, from the above inequalities, we conclude that
$$\frac{|\mathcal{E}^{(t+1)}|}{|\mathcal{E}^{(t)}|} \le \frac{\log(
n\bar{p})^{10}}{n\bar{p}} \le \frac{1}{\sqrt{n\bar{p}}}.$$
Next we prove all the steps of the previous analysis.

\noindent{\em Proof of \eqref{eq:7-1}:} From $\log(1+x)\le x$, when $p(j,i,\ell)-|\hat{p}(j,i,\ell)-p(j,i,\ell)| >0$,
$$\left|\log   \frac{\hat{p}(j,i,\ell)}{p(j,i,\ell)} \right|\le\frac{|\hat{p}(j,i,\ell)-p(j,i,\ell)|}{p(j,i,\ell)-|\hat{p}(j,i,\ell)-p(j,i,\ell)|}.$$
Thus, we just provide an upper bound of $|\hat{p}(j,i,\ell)-p(j,i,\ell)|$ to show \eqref{eq:7-1}. From the triangle inequality,
\begin{eqnarray}
&&|\hat{p}(j,i,\ell)-p(j,i,\ell)|\cr
& = & \frac{\left|e(S^{(0)}_i, S^{(0)}_j ,\ell) - p(j,i,\ell) |S^{(0)}_i| |S^{(0)}_j| \right|}{|S^{(0)}_i| |S^{(0)}_j|} \cr
&\le& \frac{\left|e(S^{(0)}_i, S^{(0)}_j ,\ell) -\mathbb{E}[e(S^{(0)}_i, S^{(0)}_j,\ell)] \right| + \left| \mathbb{E}[e(S^{(0)}_i, S^{(0)}_j,\ell)] - p(j,i,\ell) |S^{(0)}_i| |S^{(0)}_j| \right|}{|S^{(0)}_i| |S^{(0)}_j|}. \label{eq:tpp}
\end{eqnarray} 
We first find an upper bound of $\left|e(S^{(0)}_i, S^{(0)}_j ,\ell) -\mathbb{E}[e(S^{(0)}_i, S^{(0)}_j,\ell)] \right|$. Let $\mathcal{S}$ be the of partitions such that
$$\left|\cup_{k=1}^K \mathcal{V}_k \setminus S_k \right| \le \xi = O\left(\frac{\log(n\bar{p})^2}{\bar{p}} \right)\quad \mbox{for all} \quad \{S_k \}_{1\le k \le K} \in \mathcal{S}. $$
Then,
\begin{eqnarray}
\left|\mathcal{S} \right| &\le & {n \choose \xi } K^{\xi}\cr
&\le& \left(\frac{k e n}{\xi}\right)^\xi \cr
&= &\exp\left(O\left(\frac{\log(n\bar{p})^3}{\bar{p}} \right) \right).\label{eq:sizeS}
\end{eqnarray}
For all $\{S_k \}_{1\le k \le K} \in \mathcal{S}$ and for all $\ell \ge 1$ and $1\le
i,j \le K$,  $e(S_i, S_j ,\ell)$ is the sum of $|S_i||S_j|$ (or $\frac{|S_i|^2}{2}$ when $i=j$) independent Bernoulli random variables. Since the variance of $e(S_i, S_j ,\ell)$ is always less than $n^2\bar{p}$,
by Chernoff inequality (e.g., Theorem 2.1.3 in \cite{tao2012}), with probability at least $1-\exp\left(-\Theta\left(\frac{\log (n\bar{p})^4}{\bar{p}}\right) \right)$,
\begin{equation} \left|e(S_i, S_j ,\ell) - \mathbb{E}[e(S_i, S_j,\ell)]\right|  \le n \log (n\bar{p})^2 \quad\mbox{for all}\quad i,j,\ell.\label{eq:chernoffS}\end{equation}
From \eqref{eq:sizeS} and \eqref{eq:chernoffS}, with high probability,
$$\left|e(S_i, S_j ,\ell) - \mathbb{E}[e(S_i, S_j,\ell)]\right|  \le n \log (n\bar{p})^2 \quad\mbox{for all}\quad i,j,\ell \quad\mbox{and}\quad \{S_k \}_{1\le k \le K} \in \mathcal{S}. $$
Since $\{S^{(0)}_k \}_{1\le k \le K} \in \mathcal{S}$, from the above inequality,
\begin{equation} \left|e(S^{(0)}_i, S^{(0)}_j ,\ell) - \mathbb{E}[e(S^{(0)}_i, S^{(0)}_j,\ell)]\right|  \le n \log (n\bar{p})^2 \quad\mbox{for all}\quad i,j,\ell.\label{eq:6}\end{equation}

We now devote to the remaining part of \eqref{eq:tpp}. Since $|\mathcal{E}^{(0)} | =  O\left( \frac{\log (n\bar{p})^2}{\bar{p}}\right)$ from
Theorem~\ref{thm:spec}, 
\begin{equation}\left|\mathbb{E}[e(S^{(0)}_i, S^{(0)}_j,\ell)] -
  |S^{(0)}_i||S^{(0)}_j|p(i,j,\ell)\right| \le \eta |\mathcal{E}^{(0)}
|n p(i,j,\ell) =O(n \log (n\bar{p})^2).   \label{eq:7}\end{equation}
From \eqref{eq:tpp}, \eqref{eq:6} and \eqref{eq:7}, with high probability, 
$$|\hat{p}(j,i,\ell)-p(j,i,\ell)| = O(\log (n\bar{p})^2/n) \quad\mbox{for all}\quad i,j,\ell, $$
which implies that: 
\begin{eqnarray*}
\left|\log   \frac{\hat{p}(j,i,\ell)}{p(j,i,\ell)} \right|
&\le&\frac{|\hat{p}(j,i,\ell)-p(j,i,\ell)|}{p(j,i,\ell)-|\hat{p}(j,i,\ell)-p(j,i,\ell)|} = O\left(\frac{\log (n\bar{p})^2}{np(j,i,\ell)} \right) \quad\mbox{for all}\quad i,j,\ell.
\end{eqnarray*}
Since $e(v, S^{(t)}_i,\ell) \le e(v,\mathcal{V}) \le 10 \eta n\bar{p}L$ from (H1) and $np(j,i,\ell) \ge (n\bar{p})^\kappa$ from (A3), we deduce that, for all $v\in \Gamma$ and $i,j,k$,
$$\sum_{\ell=0}^L e(v, S^{(t)}_i,\ell)\left| \log  \frac{\hat{p}(j,i,\ell)}{\hat{p}(k,i,\ell)} - \log  \frac{p(j,i,\ell)}{p(k,i,\ell)} \right| = O\left(\log (n\bar{p})^2 (n \bar{p})^{1-\kappa}  \right).$$

\medskip
\noindent{\em Proof of \eqref{eq:7-2}:} Since $\log\frac{p(j,i,0)}{p(k,i,0)} = O(\bar{p})$ for all $i,j,k$ and $|\set{E}^{(t)}| = O(\log (n\bar{p})^2/\bar{p}) $,  
\begin{align*}
\sum_{i=1}^K&\sum_{\ell =0}^L  e(v, S^{(t)}_i,\ell) \log
       \frac{p(j,i,\ell)}{p(k,i,\ell)}  \cr
=&\sum_{i=1}^K \left(|S^{(t)}_i|\log\frac{p(j,i,0)}{p(k,i,0)}+ \sum_{\ell =1}^L e(v, S^{(t)}_i,\ell) \log
       \frac{p(j,i,\ell)p(k,i,0)}{p(k,i,\ell)p(j,i,0)} \right)\cr
\le&\sum_{i=1}^K \left(|\mathcal{V}_i|\log\frac{p(j,i,0)}{p(k,i,0)}+ \sum_{\ell =1}^L e(v, S^{(t)}_i,\ell) \log
   \frac{p(j,i,\ell)p(k,i,0)}{p(k,i,\ell)p(j,i,0)} \right) + \log( n \bar{p})^3\cr
\le &\sum_{i=1}^K\sum_{\ell=0}^L e(v,\mathcal{V}_i,\ell)\log\frac{p(j,i,\ell)}{p(k,i,\ell)}+ \sum_{i=1}^K\sum_{\ell
   =1}^L e(v, \mathcal{V}_i \setminus S^{(t)}_i,\ell)\log
(2\eta ) + \log( n \bar{p})^3\cr
=&\sum_{i=1}^K\sum_{\ell=0}^L
     e(v,\mathcal{V}_i,\ell)\log\frac{p(j,i,\ell)}{p(k,i,\ell)}+ \left( e(v,
     \mathcal{E}^{(t)}) + e(v,\mathcal{V}\setminus H) \right)\log
(2\eta ) + \log( n \bar{p})^3\cr
\le&\sum_{i=1}^K\sum_{\ell=0}^L
     e(v,\mathcal{V}_i,\ell)\log\frac{p(j,i,\ell)}{p(k,i,\ell)}+ \log(2\eta ) e(v,
     \mathcal{E}^{(t)}) +  2\log( n \bar{p})^3,
\end{align*}
where the last inequality stems from (H3), i.e., from $e(v,\mathcal{V}\setminus H) \le 2 \log (n\bar{p})^2$ when $v \in H$.

\medskip
\noindent{\em Proof of \eqref{eq:7-3}:} Since $\mathcal{E}^{(t+1)}
\subset H$ and every $v\in H$ satisfies (H2), every $v \in
\mathcal{E}_{jk}^{(i+1)}$ satisfies:
$$ \sum_{i =1}^K\sum_{\ell=0}^L  e(v, \mathcal{V}_i , \ell ) \log
      \frac{p(j,i,\ell)}{p(k,i,\ell)} \le -\frac{n\bar{p}}{\log (n\bar{p})^4} .$$

\medskip
\noindent{\em Proof of \eqref{eq:7-4}:} Let $\bar{\Gamma}=\{v: e(v,\mathcal{V}) \le
10\eta n\bar{p}L \}$ and $A^{\ell}_{\bar{\Gamma}}$ be the trimmed matrix of
$A^{\ell}$ whose elements in rows and columns corresponding to $w\notin \bar{\Gamma}$ are set to 0. $\bar{\Gamma}$ is the set of all items that satisfy (H1) and $H \subset \bar{\Gamma}$. Let $X_{\bar{\Gamma}} =
\sum_{\ell=1}^L (A^{\ell}_{\bar{\Gamma}}-M_{\bar{\Gamma}}^{\ell})$. We have:
$$\sum_{v \in \mathcal{E}^{(t+1)}} (e(v, \mathcal{E}^{(t)})-\mathbb{E}[e(v,\mathcal{E}^{(t)})]) \le 1_{\mathcal{E}^{(t)}}^{T} \cdot X_{\bar{\Gamma}}\cdot 1_{\mathcal{E}^{(t+1)}},
$$ 
where $1_S$ is the vector whose $v$-th component is equal to 1 if $v\in S$ and to 0 otherwise. Since $\mathbb{E}[e(v,\mathcal{E}^{(t)})]\le \bar{p} L |\mathcal{E}^{(t)}|$ and $\| X_{\bar{\Gamma}}\|_2 \le
    \sqrt{n\bar{p}\log n \bar{p}}$ with high probability from Lemma~\ref{lem:FOspectral},
\begin{align*}\sum_{v \in
      \mathcal{E}^{(t+1)}} e(v, \mathcal{E}^{(t)})=&
\sum_{v \in
      \mathcal{E}^{(t+1)}} \left(e(v, \mathcal{E}^{(t)})-\mathbb{E}[e(v,
    \mathcal{E}^{(t)})]\right) + \bar{p} L |\mathcal{E}^{(t)}||\mathcal{E}^{(t+1)}|\cr
\le & \ \|1_{\mathcal{E}^{(t)}}^{T} \cdot X_{\bar{\Gamma}}
    \cdot 1_{\mathcal{E}^{(t+1)}}\|_2 + |\mathcal{E}^{(t+1)}|\log(n\bar{p})\cr 
\le& \ \|1_{\mathcal{E}^{(t)}}^{T}\|_2 \| X_{\bar{\Gamma}}\|_2
    \| 1_{\mathcal{E}^{(t+1)}}\|_2 + |\mathcal{E}^{(t+1)}|\log(n\bar{p})\cr
\le& \ \sqrt{
      |\mathcal{E}^{(t)}||\mathcal{E}^{(t+1)}|n\bar{p} \log (n\bar{p})}+ |\mathcal{E}^{(t+1)}|\log(n\bar{p}).\end{align*}
\QED

\section{Proof of Theorem~\ref{thm:exact}}
The positive result is obtained by applying Theorem~\ref{thm:algorithms} to $s= \frac{1}{2}$. When $\lim\inf_{n\to \infty} \frac{nD(\alpha,p)}{\log(n)} \ge 1$, SP algorithm find clusters exactly with high probability.
Thus, it suffices to show the negative result. 

We prove the negative part by contradiction. Consider a maximum a posteriori (MAP) estimation with full parameter information. When we observe a labeld information $A$, the MAP estimates the clusters as follows: 
\begin{equation}
(\hat{S}_k)_{k=1,\dots,k}=\arg\max_{(S_k)_{k=1,..,K}}\mathbb{P}\left\{(S_k)_{k=1,..,K}|\alpha,p,K,A \right\}.\label{eq:MAP1}
\end{equation}
Let $\varepsilon^{\mbox{MAP}}$ denote the number of misclassified nodes by the MAP estimation.
From the definition of the MAP estimation, for any clustring algorithm $\pi$, we have
\begin{equation}
\mathbb{P}\left\{\varepsilon^\pi \ge 1 \right\} \ge \mathbb{P}\left\{\varepsilon^{\mbox{MAP}} \ge 1 \right\}.\label{eq:MAP2}
\end{equation}
Thus, in what follows, we show that when $\lim\inf_{n\to \infty} \frac{nD(\alpha,p)}{\log(n)} < 1$, the MAP estimation is failed to find the exact clusters with high probability.

We start by Lemma~\ref{lem:ex} which finds a large deviation inequality for edge connections. 
\begin{lemma}
Let $\boldsymbol{x}\in \mathbb{Z}^{K\times (L+1)}$ whose $(k,\ell+1)$ element is $x_{k,\ell}$, and such that $\sum_{\ell=0}^{L} x_{k,\ell} = |\mathcal{V}_k|$ for all $1\le k\le K$, $\sum_{\ell=1}^L x_{k,\ell}= \Theta( n \bar{p})$ for all $k$, and 
$$\sum_{k=1}^K |\mathcal{V}_k| KL(\mu(v,\mathcal{V}_k),p(i,k)) = nD \quad\mbox{when}\quad e(v) = \boldsymbol{x},$$
where we denote by $e(v)$ the $K\times
(L+1)$ matrix whose $(k,\ell+1)$ element is $e(v,\mathcal{V}_k,\ell)$. Then,
$$\log\left( \mathbb{P}\left\{ e(v) = \boldsymbol{x} \right\} \right) \ge - nD (1 + o(1))\quad\mbox{when}\quad v\in \mathcal{V}_i\quad\mbox{and}\quad D = \Omega(\bar{p}).$$  \label{lem:ex}
\end{lemma}
\noindent{\em Proof.}
When using the convention $\sum_{\ell=a}^b$ as $0$ when $a>b$, we have
\begin{align*}
& \log \left( \mathbb{P}\left\{  e(v) = \boldsymbol{x} \right\} \right)\cr
=& ~\sum_{k=1}^K\left( \left(|\mathcal{V}_k|-\sum_{\ell=1}^L x_{k,\ell}\right)\log \left(p(i,k,0) \right) + \sum_{\ell=1}^L \log\left(p(i,k,\ell)^{x_{k,\ell}}{|\mathcal{V}_k|-\sum_{m=1}^{\ell-1} x_{k,m} \choose x_{k,\ell}} \right) \right) \cr
 \ge& ~\sum_{k=1}^K\left( \left(|\mathcal{V}_k|-\sum_{\ell=1}^L x_{k,\ell}\right)\log \left(p(i,k,0) \right) + \sum_{\ell=1}^L \log\left(p(i,k,\ell)^{x_{k,\ell}}\frac{\left( |\mathcal{V}_k|-\sum_{m=1}^L x_{k,m} \right)^{x_{k,\ell}}}{ x_{k,\ell}!} \right) \right) \cr
\stackrel{(a)}{\ge}& ~\sum_{k=1}^K\left( \left(|\mathcal{V}_k|-\sum_{\ell=1}^L x_{k,\ell}\right)\log \left(p(i,k,0) \right) + \sum_{\ell=1}^L \log\left(\left( \frac{p(i,k,\ell) e}{ \frac{x_{k,\ell} }{|\mathcal{V}_k|-\sum_{m=1}^L x_{k,m}}  }\right)^{x_{k,\ell}} \frac{1}{e\sqrt{x_{k,\ell}}} \right) \right) \cr
\stackrel{(b)}{=}& ~\sum_{k=1}^K\left( \left(|\mathcal{V}_k|-\sum_{\ell=1}^L x_{k,\ell}\right)\log \left(p(i,k,0) \right) + \sum_{\ell=1}^L \log\left( \frac{p(i,k,\ell) e}{ \frac{x_{k,\ell} }{|\mathcal{V}_k|-\sum_{m=1}^L x_{k,m}}  }\right)^{x_{k,\ell}}  \right)  -o\left(\sum_{k=1}^K\sum_{\ell=1}^L x_{k,\ell}\right)\cr
\stackrel{(c)}{\ge}& ~\sum_{k=1}^K \left(|\mathcal{V}_k|-\sum_{\ell=1}^L x_{k,\ell}\right)\log  \left(p(i,k,0)\left(1+ \frac{\sum_{\ell=1}^L x_{k,\ell} }{|\mathcal{V}_k|-\sum_{\ell=1}^L x_{k,\ell}} \right) \right)  \cr
&\quad + \sum_{k=1}^K\left(\sum_{\ell=1}^L x_{k,\ell} \log\left( \frac{p(i,k,\ell) }{ \frac{x_{k,\ell} }{|\mathcal{V}_k|-\sum_{m=1}^L x_{k,m}}  }\right)  \right) -o\left(\sum_{k=1}^K\sum_{\ell=1}^L x_{k,\ell}\right)\cr
=& ~\sum_{k=1}^K \left(|\mathcal{V}_k|-\sum_{\ell=1}^L x_{k,\ell}\right)\log  \left(\frac{p(i,k,0) }{(|\mathcal{V}_k|-\sum_{\ell=1}^L x_{k,\ell})/ |\mathcal{V}_k|}  \right) + \sum_{k=1}^K\left(\sum_{\ell=1}^L x_{k,\ell} \log\left( \frac{p(i,k,\ell) }{ x_{k,\ell} /|\mathcal{V}_k|  }\right)  \right) \cr
&\quad + \sum_{k=1}^K\left(\sum_{\ell=1}^L x_{k,\ell} \log\left( \frac{ |\mathcal{V}_k|-\sum_{m=1}^L x_{k,m}  }{ |\mathcal{V}_k| }\right)  \right) -o\left(\sum_{k=1}^K\sum_{\ell=1}^L x_{k,\ell}\right)\cr
 \stackrel{(d)}{\ge}& ~ - n D -o\left(\sum_{k=1}^K\sum_{\ell=1}^L x_{k,\ell}\right) \cr
 \stackrel{(e)}{\ge}&~ - nD (1+o(1)),
\end{align*}
where $(a)$ is obtained from $n! \le e\sqrt{n} \left(\frac{n}{e} \right)^{n}$; $(b)$ stems from $\sum_{k=1}^K\sum_{\ell=1}^L x_{k,\ell} = \omega(1)$; to derive $(c)$, we use $e^{\sum_{\ell=1}^L x_{k,\ell}} \ge \left(1+ \frac{\sum_{\ell=1}^L x_{k,\ell} }{|\mathcal{V}_k|-\sum_{\ell=1}^L x_{k,\ell}} \right)^{|\mathcal{V}_k|-\sum_{\ell=1}^L x_{k,\ell}} $ since $e\ge (1+1/x)^x$ for all $x>0$; to prove $(d)$, we use the definition of ${\bm x}$ and the following inequality: 
\begin{align*}
\sum_{\ell=1}^L x_{k,\ell} \log \left(\frac{|\mathcal{V}_k|}{|\mathcal{V}_k| - \sum_{m=1}^L x_{k,m}} \right) & = \frac{\left( \sum_{\ell=1}^L x_{k,\ell} \right)^2}{|\mathcal{V}_k| - \sum_{\ell=1}^L x_{k,\ell}} (1+o(1)) = o(\sum_{\ell=1}^L x_{k,\ell});
\end{align*}
and $(e)$ is obtained from the definition of ${\bm x}$ that $\sum_{\ell=1}^L x_{k,\ell}= \Theta( n \bar{p})$ for all $k$. \QED

Assume that there exists a constant $\eta >0$ such that $\frac{nD(\alpha, p)}{\log (n)} < 1- \eta $.\\
Let $(i^\star, j^\star) = \arg\min_{i,j: i<j} D_{L+}(p(i),p(j))$ (i.e., it is the hardest case to discriminate cluster $i^\star$ and cluster $j^\star$).  When $n\to \infty$, one can easily check using the continuity of the KL divergence that there exists ${\bm x}^\star$ such that when $e(v) = \boldsymbol{x}^\star$,
\begin{align}
&\frac{\eta}{2}\log n +\sum_{k=1}^K |\mathcal{V}_k| KL(\mu(v,\mathcal{V}_k),p(j^\star,k)) < \sum_{k=1}^K |\mathcal{V}_k| KL(\mu(v,\mathcal{V}_k),p(i^\star,k)) \quad\mbox{and} \label{eq:ve-1}\\
&\sum_{k=1}^K |\mathcal{V}_k| KL(\mu(v,\mathcal{V}_k),p(i^\star,k)) \le (1-\eta/2)\log (n).\label{eq:ve-2}
\end{align}

Let $\mathcal{V}_e = \left\{v \in \mathcal{V}_{i^\star}: e(v) = \boldsymbol{x}^\star \right\}$. From \eqref{eq:ve-2} and Lemma~\ref{lem:ex}, $\mathbb{E}[|\mathcal{V}_e|] \ge n^{\eta/4}$. Thus, from Markov inequality, with probability at least $1-n^{-\eta/4}$, $\mathcal{V}_e$ is not empty (i.e., $|\mathcal{V}_e| \ge 1$). 

Let $v^\star \in \mathcal{V}_e $ be a node in $\mathcal{V}_e$.
We denote by $\Phi$ the original partition and define a slightly modified partition $\Psi$ as follows:
$$\hat{\mathcal{V}}_{i^\star} = \mathcal{V}_{i^\star} \setminus \{ i^\star \}, \quad \hat{\mathcal{V}}_{j^\star} = \mathcal{V}_{j^\star} \cup \{ i^\star \},\quad\mbox{and}\quad
\hat{\mathcal{V}}_k = \mathcal{V}_k\quad\mbox{otherwise}.$$
Then, $\Psi$ is a more likely partition than $\Phi$ from \eqref{eq:ve-1}, i.e.,
\begin{equation}
\mathbb{P}\left\{\Phi|\alpha,p,K,A \right\} \ge \mathbb{P}\left\{\Psi|\alpha,p,K,A \right\}
\end{equation}
which means that the MAP estimator does not select the exact partition when $\mathcal{V}_e$ is not empty. Therefore, from \eqref{eq:MAP2}, every clustering algorithm $\pi$ has the error probability that
$$\mathbb{E}\left\{\varepsilon^\pi \ge 1 \right\} \ge 1-n^{-\eta/4} $$
when there exists a constant $\eta >0$ such that $\frac{nD(\alpha, p)}{\log (n)} < 1- \eta $. \QED

\end{document}